\documentclass{article}

\usepackage{arxiv}

\usepackage[utf8]{inputenc} 
\usepackage[T1]{fontenc}    

\usepackage{booktabs}       
\usepackage{threeparttable}
\usepackage{tabularx}

\usepackage{graphicx}
\usepackage{tikz} 
\usetikzlibrary{arrows.meta, positioning}

\usepackage{amsfonts}       
\usepackage{amssymb}
\usepackage{amsmath}
\usepackage{bm} 
\usepackage{nicefrac}       

\usepackage[
	locale=US,
	per-mode = fraction,
	range-phrase = \; to \;,
	range-units = repeat
	]{siunitx}  

\usepackage{algorithm}
\usepackage[noend]{algpseudocode}

\usepackage[printonlyused]{acronym}

\usepackage[%
    colorlinks=true,
    pdfborder={0 0 0},
    citecolor=blue,
    linkcolor=blue,
    urlcolor=blue
]{hyperref}
\usepackage[nameinlink]{cleveref} 
\usepackage{url}            

\usepackage{microtype}      
\usepackage[square,numbers]{natbib}
\usepackage{doi}

\input{macros.sty}

\title{Non-intrusive hyperreduction by a physics-augmented neural network with second-order Sobolev training}

\author{ \href{https://orcid.org/0000-0003-3285-7386}{\includegraphics[scale=0.06]{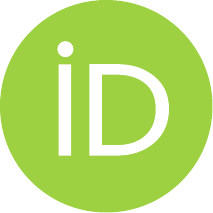}\hspace{1mm}Arwed Schütz} \\
	Jade University of Applied Sciences\\
	Wilhelmshaven, Germany \\
	\texttt{arwed.schuetz[at]jade-hs.de} \\
	\And
	\href{https://orcid.org/0009-0000-8796-7764}{\includegraphics[scale=0.06]{orcid.pdf}\hspace{1mm} Lars Nolle} \\
	Jade University of Applied Sciences\\
	Wilhelmshaven, Germany \\
    German Research Centre for Artificial Intelligence\\
	\texttt{lars.nolle[at]jade-hs.de} \\
	\And
	\href{https://orcid.org/0000-0003-1993-4893}{\includegraphics[scale=0.06]{orcid.pdf}\hspace{1mm} Tamara Bechtold} \\
	Jade University of Applied Sciences\\
	Wilhelmshaven, Germany \\
	\texttt{tamara.bechtold[at]jade-hs.de} \\
}

\renewcommand{\headeright}{}
\renewcommand{\shorttitle}{}

\hypersetup{
pdftitle={Non-intrusive hyperreduction by a physics-augmented neural network with second-order Sobolev training},
pdfauthor={Arwed Schütz, Lars Nolle, Tamara Bechtold}
}

\begin{document}
\maketitle

\begin{abstract}
The finite element method is an indispensable tool in engineering, but its computational complexity prevents applications for control or at system-level. Model order reduction bridges this gap, creating highly efficient yet accurate surrogate models.  Reducing nonlinear setups additionally requires hyperreduction. Compatibility with commercial finite element software requires non-intrusive methods based on data. Methods include the trajectory piecewise linear approach, or regression, typically via neural networks. Important aspects for these methods are accuracy, efficiency, generalization, including desired physical and mathematical properties, and extrapolation. Especially the last two aspects are problematic for neural networks. Therefore, several studies investigated how to incorporate physical knowledge or desirable properties. A promising approach from constitutive modeling is physics augmented neural networks. This concept has been elegantly transferred to hyperreduction by Fleres et al. in 2025 and guarantees several desired properties, incorporates physics, can include parameters, and results in smaller architectures. We augment this reference work by second-order Sobolev training, i.e., using a function and its first two derivatives. These are conveniently accessible and promise improved performance. Further modifications are proposed and studied. While Sobolev training does not meet expectations, several minor changes improve accuracy by up to an order of magnitude. Eventually, our best model is compared to reference work and the trajectory piecewise linear approach. The comparison relies on the same numerical case study as the reference work and additionally emphasizes extrapolation due to its critical role in typical applications. Our results indicate quick divergence of physics-augmented neural networks for extrapolation, preventing its deployment.
\end{abstract}


\clearpage

\section{Introduction}

Real-world physics ranging from structural mechanics to electromagnetism are described by \acp{PDE}.
However, analytical solutions are limited to trivial scenarios that are hardly related to modern technical products.
A well-established method to bridge this gap is the \ac{FEM}.
It allows for performing design studies, optimization, and reliability analysis before fabrication to minimize the number of physical validations.
Furthermore, all quantities of interest are available at every position and sampling frequency, whereas experimental measurements are limited by sensors, environmental conditions, and accessibility.
Although \ac{FEM} is a powerful tool, underlying models are high-dimensional and associated with high computational costs.
This computational burden prevents \ac{FEM} from several applications~\cite{Hartmann2018}, e.g., in system-level simulations or feedback control.

An established method to overcome this challenge for linear scenarios is projection-based \ac{MOR}.
Based on the original model, \ac{MOR} generates a surrogate model of drastically smaller dimension while hardly compromising on accuracy~\cite{Antoulas.2005}.
In an intuitive sense, these techniques identify patterns in the solution.
The original model is then reduced by approximating the solution as a linear combination of these patterns and additionally projecting the system onto a suitable subspace.
Therefore, \ac{MOR} offers physical intuition as all results can be transformed back to the original model.
A prerequisite is access to the underlying system matrices, which is generally fulfilled in the context of \ac{FEM}, even for commercial software~\cite{mor4ansys}.
Note that methods without this prerequisite exist~\cite{Peherstorfer.2016,Proctor2016, Kneifl2021}, but for the context of \ac{FEM}, the reference quantities are directly available.

However, \ac{MOR} faces a significant computational bottleneck when applied to nonlinear systems.
In order to evaluate nonlinear terms, the reduced state has to be expanded back to its original dimension, then the same set of nonlinearities as in the original model must be evaluated, which are then projected back to the reduced dimension.
To resolve this bottleneck, \ac{MOR} is augmented by additional methods referred to as \textit{hyperreduction}~\cite{Ryckelynck2005b} which efficiently approximate nonlinearities at reduced level~\cite{Baur.2014}. 
Prominent hyperreduction methods are the \ac{DEIM}~\cite{Chaturanta2010} or the \ac{ECSW} method~\cite{Farhat2014}.
Both methods rely on evaluating a small subset of the original nonlinear terms.
Therefore, they require deep access to the original model's code and are referred to as \textit{intrusive} methods.
As this level of access is hardly available in most \ac{FEM} software, \textit{non-intrusive} methods have gained attention.

Non-intrusive hyperreduction has a long history, including curve fitting or interpolating linearizations, but also modern approaches based on \acp{ANN}.
This work classifies all methods as non-intrusive that exclusively rely on data that can be extracted from commercial \ac{FEM} software such as Abaqus~\cite{Hollkamp2008}, Ansys~\cite{mor4ansys}, COMSOL~\cite{Nokhosteen2019}, Nastran~\cite{Phlipot2014}, etc.
A method established early is energy-based hyperreduction~\cite{Gabbay1998}.
This approach samples energies for different reduced states and fits a rational polynomial.
Its gradient then constitutes the reduced forces.
Furthermore, its Hessian resembles the tangent stiffness matrix which is symmetric if the second partial derivatives are continuous. 
This procedure has been successfully applied to micro-electromechanical systems~\cite{Gabbay1998,Gabbay.2000,Mehner.2000,Bennini2001} or magnetic devices~\cite{Varghese1999}.
Another early and simple method is the \ac{TPWL} approximation \cite{MichalJerzyRewienski.2003}.
The strategy is to dynamically interpolate linearizations.
The weighting scheme relies on the distances between the current state and the states corresponding to the linearizations.
\ac{TPWL} is robust and relies only on samples of state vector, force vector, and the tangent stiffness matrix.
These data are readily available in the context of \ac{FEM} as it solves nonlinear problems via iterative linearizations.
A modern approach for non-intrusive hyperreduction are \acp{ANN}~\cite{Gao2020}.
In its basic form, this approach is a regression problem to approximate the relation between reduced state and reduced forces.
However, this method constitutes a black box that maps data without physical reasoning or guaranteed properties.
Note that polynomial alternatives exist~\cite{Boef2024}, but these do not enjoy the universal approximation theorem~\cite{Hornik1989}.

Recently, a highly promising \ac{ANN} method that mitigates these drawbacks has been transferred from constitutive modeling~\cite{Asad2022,Klein2022,Linden2023} to hyperreduction~\cite{Fleres2025} and augmented to also capture parametric effects.
Similarly to early energy-based hyperreduction~\cite{Gabbay1998}, the \ac{ANN} maps the reduced state to an energy-like quantity.
The function's gradient provides the reduced forces on which the \ac{ANN} is trained.
This concept is related to Lagrangian~\cite{Cranmer2020} or Hamiltonian neural networks~\cite{Greydanus2019} or work on potential energy surfaces for molecular dynamics~\cite{Blank1995,Behler2011}.
Another advantage is a symmetric Hessian.
In addition, deploying an \ac{ICNN}~\cite{Amos2017} ensures a convex function, guaranteeing a positive semi-definite Hessian and thus, a stable \ac{ROM}.
This is a notable improvement over other work that only encourages stability as an additional loss term~\cite{Sharma2024}.
Augmenting the \ac{ANN} by a linear offset term ensures mechanical consistency~\cite{Asad2022}, i.e., zero forces for zero displacement.
Hence, this type of approach elegantly embeds a priori knowledge and guarantees many desirable properties by design.
In contrast to the well-established \acp{PINN}~\cite{Raissi2017,Kharazmi2019,Haghighat2021}, which penalize violations of physics during training, physics or properties are not promoted but guaranteed.
Corresponding implementations are also referred to as \acp{PANN}~\cite{Klein2022b}.

In summary, the recent \ac{PANN}-based hyperreduction method~\cite{Fleres2025} has many advantages:
Its non-intrusive nature ensures compatibility even with commercial \ac{FEM} software.
The physics-augmented approach does not suggest but guarantees desired properties and improves generalization.
In addition, the concept of energy elegantly introduces physical interpretation, also for derivatives as force and tangent stiffness.
Using scalar-valued energy as output instead of vector-valued force also results in a smaller network, reducing computational costs of training and evaluation.
Moreover, parametric capabilities can be included, e.g., for material parameters.
Successful related work in other branches emphasizes the concept's potential.

However, the approaches so far neglect information commonly present in the context of nonlinear \ac{FEM}: energies and tangent stiffness matrices.
While the former are available as postprocessing quantities, the latter arise as byproducts of the \ac{FEM} solution process via iterative linearizations.
Including derivatives in regression can enhance accuracy, data-efficiency, and generalization~\cite{Czarnecki2017}.
For these reasons, the general idea has a long history with several methods, including regression techniques~\cite{Laurent2017} or \acp{ANN}. 
In the field of regression, surrogate modeling, or metamodeling, corresponding methods are typically known as \textit{gradient-enhanced}~\cite{Laurent2017} methods.
Examples include Kriging~\cite{Chung2002,Han2013,Bouhlel2019}, support vector machines~\cite{Jayadeva2006}, or radial basis function networks~\cite{Kampolis2004,Giannakogl2006}.
In the more specific area of \acp{ANN}, research considering derivatives started as early as 1990~\cite{Hornik1990}.
Subsequent work~\cite{Simard1991,Drucker1992} utilized derivatives in training to regularize or to match target derivatives~\cite{Flake1999}.
Another early work~\cite{Dugas2000} introduced the \textit{Softplus} activation function and constrained weights to ensure derivative-based characteristics of the original function.
These findings are highly relevant for \acp{ICNN}~\cite{Amos2017}, an important element in \ac{ANN}-based constitutive modeling and hyperreduction.
Matching actual target derivatives has been deployed to approximate an aircraft’s trim map~\cite{Ferrari2005} or to learn potential energy surfaces~\cite{Witkoskie2005} in chemistry.
Later work investigated how to combine the losses and proposed pruning to avoid certain overfitting effects~\cite{Pukrittayakamee2011}.
An important publication in this area~\cite{Czarnecki2017} introduced the term Sobolev training and also proposed a more efficient stochastic version.
The exponentially growing data dimension limits the derivative order, but some studies have incorporated derivatives up to the fifth order~\cite{Avrutskiy2018,Avrutskiy2020}.
Another aspect is how to combine the multiple losses with respect to different goals, e.g., by scheduling, weighting, or tuning the amount of data.
While some studies suggest to gradually include derivative information during training to enhance efficiency~\cite{Bouhlel2020}, other studies~\cite{Kissel2020} propose the opposite for more accurate output predictions.
Other work studies the effect of not using all but only a fraction of samples~\cite{Feng2022c}.
These studies conclude that additionally considering derivative information improves accuracy, data-efficiency, generalization, and robustness against noise at the cost of longer training.
Work using separate outputs or even \acp{ANN}~\cite{Cicci2022} to approximate function and derivatives does not enjoy these benefits.
Consequently, the concept of Sobolev training has diffused into physically relevant applications, including \acp{PINN}~\cite{Yu2022} or a discrete time dynamical system of a robot arm~\cite{Kim2022}.
To our best knowledge, its combination with \ac{ANN}-based hyperreduction was first proposed in our earlier work~\cite{Schuetz2024b}.
However, note that Sobolev training is not guaranteed to improve performance~\cite{Schommartz2025}.

This article studies several enhancements such as Sobolev training for \ac{PANN}-based hyperreduction and compares them with the reference work~\cite{Fleres2025} and \ac{TPWL}.
For comparability, the reference numerical case study of a nonlinear static cantilever beam is reproduced with the same load case.
However, we discard parametric influences to focus solely on nonlinearities.
Moreover, additional load cases assess extrapolation capabilities.
These play a crucial role for typical applications such as system-level simulation or control, as models likely encounter states that were not included in training.
In contrast to the reference work~\cite{Fleres2025}, this work evaluates each \ac{ANN} variant using multiple initializations to account for the stochastic nature of \acp{ANN}.
Consequently, hyperparameter studies are repeated and extended to identify the most accurate model.
The overall goal of the proposed enhancements is to include as much physical knowledge as possible and to implement best practices of machine learning.
Specific enhancements include standardizing inputs, physical initializations of the quadratic pass-through layer, an additional offset to ensure consistent energy, and second-order Sobolev training.
Second-order Sobolev training matches the function, its gradient, and its Hessian.
In a physical sense, these correspond to energy, force, and tangent stiffness, which are conveniently accessible even with commercial \ac{FEM} software~\cite{Phlipot2014,Hollkamp2008,mor4ansys,Nokhosteen2019,Perez2014}.
The potential benefits of improved accuracy, data-efficiency, and generalization are important in the context of \ac{FEM} where samples are typically scarce compared to the model's dimension.
Note that Sobolev training introduces additional loss terms, which potentially must be balanced to compensate for different orders of magnitude or convergence rates.
Therefore, various approaches are examined for combining the losses effectively.

This paper is structured as follows:
\Cref{sec:methods} describes the methodology, including \ac{FEM}, \ac{MOR}, \ac{TPWL}, and the \ac{PANN} and its aspects.
\Cref{ssec:Results} applies the methodology to a numerical case study and presents the corresponding results.
The numerical case study is a static nonlinear cantilever beam reproduced from the reference work~\cite{Fleres2025}.
After reporting on the \ac{FEM} model, data collection, and dimensionality reduction, the \ac{ANN}-based hyperreduction is investigated in two stages: in the first stage, the \ac{ANN} is trained following the reference work and the best set of hyperparameters identified.
The second stage adds second-order Sobolev training and investigates multiple strategies to optimally balance the three losses for energy, force, and stiffness.
The best performing \ac{ANN} is compared with the results of reference work~\cite{Fleres2025} and the \ac{TPWL} model.
In contrast to the reference work, the comparison additionally investigates extrapolation because of its importance for typical applications such as control or system-level simulation.
\Cref{sec:Conclusion} summarizes our findings and \Cref{sec:Outlook} proposes future steps.

All scripts are publicly available at GitLab at \url{https://gitlab.gwdg.de/jade-hochschule/fms/2025_cmame}, including the \ac{FEM} model, data dimensionality reduction, \ac{PANN} training, visualization of the hyperparameter studies, and the comparison of hyperreduced models.

\section{Methodology}
\label{sec:methods}

The methodology covers the whole modeling workflow, ranging from \ac{FEM} over \ac{MOR} to hyperreduction via \ac{TPWL} or a designated \ac{ANN}.
A brief summary of nonlinear static \ac{FEM} is presented in \Cref{ssec:FEM}.
Subsequently, \Cref{ssec:MOR} describes how to reduce such a model via the well-established \ac{POD}.
The last step is non-intrusive hyperreduction, for which two approaches are detailed: \ac{TPWL} in \Cref{ssec:TPWL} and \ac{PANN} in \Cref{ssec:PANN}, including settings, architecture, and aspects of Sobolev training.
Note that the discussion focuses on static systems but is not limited to them.

\subsection{Nonlinear structural finite element models}
\label{ssec:FEM}

\ac{FEM} is a numerical method to solve differential equations based on spatial discretization.
For a static nonlinear structural case, \ac{FEM} results in the following large-scale system of nonlinear algebraic equations
\begin{equation}
\Sigma =
\begin{cases}
\bm{f} \left( \bm{x} \right) = \bm{B} \, \bm{u} \\
\bm{y} = \bm{C} \, \bm{x}
\end{cases}.
\label{EQ_original_system}
\end{equation}
In the case of structural problems, the first equation corresponds to a force balance and the second contains postprocessing.
The nonlinear restoring forces are denoted as $\vi{f}\left(\vi{x}\right)\mydim{n}$ and depend on the state vector $\vi{x}\mydim{n}$ that contains all nodal displacements.
In general, the exact function $\vi{f}\left(\vi{x}\right)\mydim{n}$ is not accessible to the user, especially in case of commercial software.
The product of the input matrix $\mi{B}\mydims{n}{p}$ and the input vector $\vi{u}$ gives $p$ external loads.
The former contains the loads' spatial distributions, the latter their magnitudes.
The $q$ user-defined outputs $\vi{y}\mydim{q}$ are computed from the state vector by the output matrix $\mi{C}\mydims{q}{n}$.

Note that the nonlinear equations as in \Cref{EQ_original_system} are typically solved via the Newton-Raphson procedure, i.e., by iterative linearizations.
Consequently, the tangent system matrices generated in the process are available to the user.
This also holds true for commercial software via features such as linear perturbation.
Along with states, forces, and energies, this information provides the sole means of accessing the nonlinear behavior within industrial software.

\subsection{Projection-based model order reduction}
\label{ssec:MOR}

\ac{MOR} constructs surrogate models of drastically smaller dimension while hardly sacrificing accuracy.
Resulting \acp{ROM} conquer applications such as feedback control or system-level simulation that were previously prevented by computational constraints.
The key idea of \ac{MOR} is that the original state vector can be approximated by a low-dimensional subspace.
In a physical sense, this subspace is spanned by dominant patterns that characterize the state's evolution, e.g., eigenmodes in case of linear structural dynamics.
Collecting these patterns as columns of a matrix and orthonormalizing them results in a projection matrix $\mi{V} \mydims{n}{r}$.
Subsequently, the original state $\vi{x}$ is approximated as
\begin{equation}
    \vi{x} \approx \mi{V} \, \vi{x}_r,
\label{EQ_MOR_state_approximation}
\end{equation}
where $\vi{x}_r\mydim{r}$ is the reduced state vector.
Hence, its components are generalized coordinates that contain the amplitudes that scale the dominant patterns in $\mi{V}$.

Especially for linear models, there are multiple methods to identify the reduced basis.
However, these typically rely on control-theoretic concepts that are not available for nonlinear setups.
A highly popular alternative, particularly for nonlinear setups, is \ac{POD}, which extracts a reduced basis from trajectories of the original model.
In detail, the original model is analyzed for representative load cases.
The resulting state vectors for all solution steps are collected as columns of a snapshot matrix.
A \ac{SVD} of this matrix then reveals dominant spatial patterns in the form of left singular vectors.
These are sorted in descending order by their importance which is indicated by their singular values.
Typically, the singular values are scaled to a sum of one for a general scale.
Note that a \ac{POD} basis depends on the data it was extracted from.
These preparations, such as data collection, training, and computing operators, are termed the \textit{offline} phase.
The subsequent \textit{online} phase refers to the model's deployment.

Once the projection matrix has been constructed, the original model's dimension is reduced in two steps: inserting \Cref{EQ_MOR_state_approximation} into \Cref{EQ_original_system} and projecting it orthogonally onto the same subspace.
The resulting reduced system reads
\begin{equation}
    \Sigma_r =
    \begin{cases}
        \mi{V}^\top \, \vi{f}\left(\mi{V} \, \vi{x}_r\right) 
        = \overbrace{\mi{V}^\top \, \mi{B}}^{\mi{B}_r} \, \vi{u} \\
        \vi{y} = \underbrace{\mi{C} \, \mi{V}}_{\mi{C}_r} \, \vi{x}_r
    \end{cases}.
\label{EQ_reduced_system}
\end{equation}
The subscript $r$ indicates reduced quantities of dimension $r\ll n$, which result from multiplication with the projection matrix.
Note that inputs and outputs remain unchanged.
However, a significant computational bottleneck remains due to the nonlinearity as this term cannot be reduced.
Instead, its evaluation requires expanding the reduced state back to the original dimension, evaluating the full expression, and projecting it back to the reduced dimension.
Therefore, evaluating nonlinear terms for the reduced model in \Cref{EQ_reduced_system} is more complex than for the original model.
Furthermore, a nonlinear model is solved via iterative linearizations, each of which requires linearizing the model.

Methods to mitigate this problem are referred to as hyperreduction.
The key idea is to find an efficient approximation or the relation between reduced forces and reduced state, i.e.,
\begin{equation}
    \mi{V}^\top \, \vi{f}\left(\mi{V} \, \vi{x}_r\right) \approx 
    \vi{f}_r \left(\vi{x}_r\right) \, .
\label{eq:hyperreduction}
\end{equation}
Well-established techniques include \ac{DEIM} and \ac{ECSW}.
Conceptually, \ac{DEIM} approximates the nonlinearities in the original dimension and subsequently reduces them while \ac{ECSW} directly approximates the reduced quantities.
Therefore, the approaches are classified as \textit{project-then-approximate} and \textit{ap\-prox\-i\-mate-then-project} methods, respectively.
Both are based on evaluating the same nonlinearity as in the original model but only a drastically smaller subset.
Hence, both techniques are classified as intrusive.
However, the corresponding expression is rarely accessible, especially in the case of commercial software.
Non-intrusive methods based on data offer an alternative, examples of which are described in the following subsections.

\subsection{Hyperreduction via trajectory-piecewise linear approximation}
\label{ssec:TPWL}

\ac{TPWL} is a non-intrusive hyperreduction method that dynamically interpolates linearizations~\cite{MichalJerzyRewienski.2003}.
Similar to \ac{POD}, these are obtained by sampling representative trajectories.
The linearizations are assembled into an effective linear model as a weighted sum.
The weights depend on the distance between the current state and the ones corresponding to the samples.
A softmax-like function scales the weights into the appropriate range and ensures sharp transitions.
This weighting scheme resembles a smooth nearest neighbor algorithm and exhibits similarities to inverse distance weighting.
\ac{TPWL} is robust and straightforward to implement, and is therefore chosen as a baseline.
Furthermore, it synergizes well with gain scheduling for designing feedback control~\cite{Tonkens2021,Schuetz2025}.

From a mathematical perspective, the first step is to linearize the nonlinear term $\vi{f}\left(\vi{x}\right)$ around $N$ expansion points $\vi{x}_i$ as
\begin{equation}
    \begin{aligned}
    \vi{f}\left(\vi{x}\right)\big|_{\vi{x}_i} 
    & \approx \underbrace{\vi{f}\left(\vi{x}\right) \big|_{\vi{x}_i}}_{\vi{f}_i}
    && + \underbrace{\dfrac{\partial \vi{f}\left(\vi{x}\right)}{\partial \vi{x}} \Big|_{\vi{x}_i}}_{\mi{K}_i} \; \left(\vi{x} - \vi{x}_i \right) \\
    & = \underbrace{\vi{f}_i - \mi{K}_i \; \vi{x}_i}_{\tilde{\vi{f}}_i}
    && + \mi{K}_i \; \vi{x}, \\
    \end{aligned}
    \label{EQ_Linearization}
\end{equation}
where $\vi{f}_i$ and $\mi{K}_i$ are the function and Jacobian evaluated at the expansion point, respectively.
The former constitutes a force vector, the latter a tangent stiffness matrix.
During operation, the $N$ linearizations are dynamically combined as a weighted sum.
The global approximation reads
\begin{equation}
\vi{f} \left( \vi{x} \right) 
\approx \sum_{i=1}^N w_i(\vi{x}) \; \tilde{\vi{f}}_i + \sum_{i=1}^N w_i(\vi{x}) \; \mi{K}_i  \; \vi{x} \, .
\label{EQ_Linearization_Combination}
\end{equation}
This approximation can be inserted into the nonlinear model and all its components are compatible with \ac{MOR}.
The weights $w_i(\vi{x})$ are state-dependent and computed according to \Cref{ALG_TPWL_weighting}.
Firstly, the distance of the current state to the linearizations' expansion points is determined.
Secondly, a softmax-like function obtains weights in the desired range.
A numerical offset $\varepsilon$ avoids singularity, and a parameter $\beta$ determines the sharpness of the transitions.

\begin{center}
    \begin{minipage}{0.6\linewidth}
    \begin{algorithm}[H]
    \caption{Weighting scheme for \ac{TPWL}.} 
    \label{ALG_TPWL_weighting}
    \begin{algorithmic}
    \For{$i = 1, \dots, N$}
        \State $d_i = \norm{\vi{x}_r - \vi{x}_{r,i} }$
    \EndFor
    \State $m = (\textrm{min}_{i=1,..,N} \; d_i) + \varepsilon$
    %
    \For{$i = 1, \dots, N$}
        \State $\hat{w_i} = \exp({-\beta \, d_i / m})$
    \EndFor
    \State $S = \sum_i^N \hat{w_i}$
    %
    \For{$i = 1, \dots, N$}
        \State $w_i = \hat{w_i} / S$
    \EndFor
    \end{algorithmic}
    \end{algorithm}
    \end{minipage}
\end{center}

For subsequent evaluations, data usage must be comparable with \acp{ANN}.
For consistency reasons, the same percentage of training data is used, but with equidistant samples from the union of the train and the validation sets.

\subsection{Hyperreduction via a physics-augmented neural network}
\label{ssec:PANN}

This work utilizes the \ac{POD}-\ac{ICNN} framework~\cite{Fleres2025} for \ac{ANN}-based hyperreduction due to its numerous advantages, but discards its parametric capabilities.
Note that the non-parametric \ac{ANN} design has been originally proposed for constitutive modeling~\cite{Asad2022,Klein2022}.
This design combines several remarkable ideas to enhance the physical interpretation and guarantee desired physical or mathematical properties.
The concept is termed \ac{PANN}~\cite{Klein2022b} to differentiate itself from \acp{PINN} as these only implement soft constraints.
The first idea is a strain energy-like scalar quantity as the \ac{ANN}'s output.
As a result, its gradients resemble restoring forces, and its Hessians tangent stiffness matrices.
Furthermore, the latter are symmetric by design due to the symmetry of second derivatives.
This idea also appears in Lagrangian~\cite{Cranmer2020} and Hamiltonian~\cite{Greydanus2019} \acp{ANN}, \acp{ANN} to model potential energy surfaces for molecular dynamics~\cite{Blank1995,Behler2011}, and also for hyperreduction~\cite{Sharma2024}.
The second idea is to use the \ac{ICNN}~\cite{Amos2017} architecture.
This choice guarantees the \ac{ANN} to be a convex function, which ensures positive semi-definite Hessians and thus numerical stability.
The third idea is a linear correction term subtracted from the \ac{ICNN}'s output to guarantee consistency, i.e., zero forces for zero displacement.
The linear offset's coefficients correspond to the \ac{ICNN}'s gradient evaluated at zero input.
In summary, these ideas elegantly incorporate a priori knowledge and improve generalization.
This is particularly relevant in the context of \ac{FEM}, where samples are typically scarce and originate from only a few load cases.
Furthermore, not all possible laodcases and disturbances faced during operation can be anticipated.

The remainder of this subsection provides details on the specific \ac{ANN} design and its components.
Furthermore, modifications are proposed, including a preceding frozen layer that standardizes the inputs, an additional constant offset to ensure consistent energy, and initializing the quadratic pass-through layer as the tangent stiffness matrix at rest.
In addition, Sobolev training, i.e., learning a function and its derivatives, and techniques to balance multiple losses are investigated.

\subsubsection{Architecture}

The basic architecture is an \ac{ICNN}~\cite{Amos2017}, which is a scalar-valued convex function.
To derive conditions for convexity, the \ac{ANN} is interpreted as a composition of individual functions, i.e., linear layers and activation functions.
A composition is convex if the initial function is convex and subsequent functions are convex and non-decreasing~\cite[p.~86]{Boyd2004}.
As a result, all linear layers but the first must have non-negative weights and all activation functions must be convex and non-decreasing.
To enhance representation power, linear pass-through layers or skip-connections~\cite{Hoedt2023} are introduced, each directly connecting the input to hidden layers.
In addition, any convex function can be added in parallel, directly connecting the input to the output.
Previous work~\cite{Asad2022,Fleres2025} has chosen a quadratic function with matrix-valued weights, inspired by energy computation in linear models.

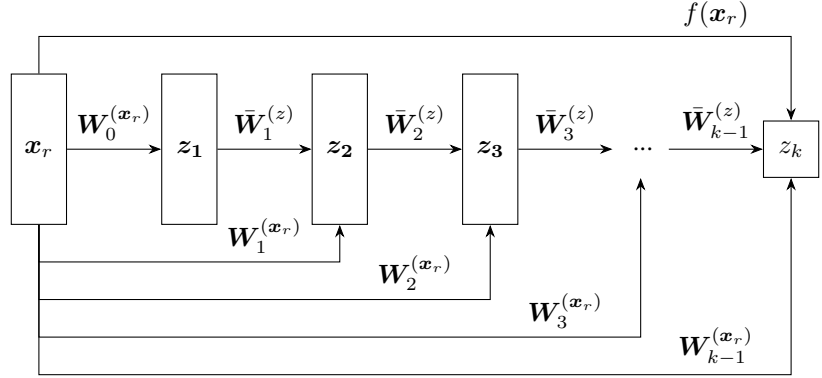
\begin{figure}[!htb]%
    \centering
    \begin{tikzpicture}[node distance=2cm, on grid, auto]
    \tikzstyle{state} = [rectangle, draw, text width=0.5cm, align=center, minimum height=2cm]
    \tikzstyle{endstate} = [rectangle, draw, text width=0.5cm, align=center, minimum height=0.75cm]
    \tikzstyle{arrow} = [draw, -{Stealth}]

    \node[state] (A) {$\vi{x}_r$};
    \node[state, right=of A] (B) {$\vi{z_1}$};
    \node[state, right=of B] (C) {$\vi{z_2}$};
    \node[state, right=of C] (D) {$\vi{z_3}$};
    \node[endstate, right=of D, draw=none] (E) {...};
    \node[endstate, right=of E] (F) {$z_k$};

    \draw[arrow] (A) -- node[midway, above]{$\mi{W}_0^{(\vi{x}_r)}$} (B);
    \draw[arrow] (B) -- node[midway, above]{$\bar{\mi{W}}_1^{(z)}$} (C);
    \draw[arrow] (C) -- node[midway, above]{$\bar{\mi{W}}_2^{(z)}$} (D);
    \draw[arrow] (D) -- node[midway, above]{$\bar{\mi{W}}_3^{(z)}$} (E);
    \draw[arrow] (E) -- node[midway, above]{$\bar{\mi{W}}_{k-1}^{(z)}$} (F);

    \draw[arrow] (A) |- ++(0,-1.5) -| node[above,xshift=-1cm]{$\mi{W}_1^{(\vi{x}_r)}$} (C);
    \draw[arrow] (A) |- ++(0,-2) -| node[above,xshift=-1cm]{$\mi{W}_2^{(\vi{x}_r)}$} (D);
    \draw[arrow] (A) |- ++(0,-2.5) -| node[above,xshift=-1cm]{$\mi{W}_3^{(\vi{x}_r)}$} (E);
    \draw[arrow] (A) |- ++(0,-3) -| node[above,xshift=-1cm]{$\mi{W}_{k-1}^{(\vi{x}_r)}$} (F);
    \draw[arrow] (A) |- ++(0,1.5) -| node[above,xshift=-1cm]{$f(\vi{x}_r)$} (F);
\end{tikzpicture}
    \caption{\ac{ICNN} architecture. Quantities marked with a bar are constrained to be non-negative. The $f(\vi{x}_r)$-term represents a quadratic energy layer, i.e., energy as it would be computed for a linear system. Note that bias terms are not shown and are not subject to constraints. The final output $z_k$ correspond to an energy-like quantity $\hat{e}$.}
    \label{fig:ICNN_architecture}
\end{figure}

\Cref{fig:ICNN_architecture} presents the \ac{ICNN} architecture and \Cref{alg:ICNN_forward} details the computation.
The input $\vi{x}_r$ propagates through the layers and results in the scalar output $z_k = \hat{e}$.
All linear layers are denoted as $\vi{W}$.
A bar denotes weights that are constrained to be non-negative.
Bias terms $\vi{b}$ and activation functions $\vi{\phi}$ are not shown, but included in \Cref{alg:ICNN_forward}.
The function $f(\vi{x}_r)$ represents possible convex functions such as the quadratic layer.

\begin{center}
    \begin{minipage}{0.6\linewidth}
    \begin{algorithm}[H]
    \caption{Forward pass.} 
    \label{alg:ICNN_forward}
    \begin{algorithmic}[1]
    \Require $\vi{x}_r$
    \Ensure $z_k = \hat{e}$
    \State $z_1 = \phi \left( \mi{W}_0^{y} \, \vi{x}_r + b_0 \right)$
    \For{$i = 1, \dots, k-1$}
        \State $\vi{z}_{i+1} = \phi \left( \bar{\mi{W}}_i^{z} \, \vi{z}_i + \mi{W}_i^{y} \, \vi{x}_r + \vi{b}_i \right)$
    \EndFor

    \State $z_k = z_k + f(\vi{x}_r)$
    
    \end{algorithmic}
    \end{algorithm}
    \end{minipage}
\end{center}

\subsubsection{Activation function}

To guarantee the \ac{ANN}'s convexity, its activation functions $\phi$ must be convex and non-decreasing~\cite{Amos2017}. 
Furthermore, they must be at least twice differentiable with non-vanishing second derivatives~\cite{Fleres2025} to compute gradients and Hessians.
Few activation functions meet these requirements, including Softplus~\cite{Dugas2000}, SoftplusSquared~\cite{Asad2022}, and SoftExponential~\cite{Godfrey2015}.
Following the reference work, this study deploys the SoftplusSquared function which is based on the Softplus function.

The Softplus function is a smooth version of the popular \ac{ReLU}.
There is a historical link as it was originally proposed to include a priori knowledge on derivatives~\cite{Dugas2000} , i.e., convexity and non-decreasing behavior.
Later work~\cite{Glorot2011} has introduced a sharpness parameter $\alpha$, resulting in its modern form 
\begin{equation}
    \text{Softplus}(x) = \dfrac{1}{\alpha} \, \text{ln} \left( 1+\exp(\alpha \, x)\right) \, .
    \label{eq:Softplus}
\end{equation}
Note that $\alpha$ must be positive in order for the function to be convex and defined.
One approach to ensure non-negative values is to square the parameter, i.e., $\alpha = \beta^2$~\cite{Asad2022}, although alternatives are available.
A drawback of Softplus is its second derivative's limited representation power as it resembles a sharp Gaussian.
This limitation may impede learning nonlinear Hessians.
Deploying large layer widths may provide a remedy in the sense of a point-wise approximation, but at the cost of a significantly larger architecture.
An alternative to overcome that shortcoming is the SoftplusSquared~\cite{Asad2022} activation function
\begin{equation}
    \text{SoftplusSquared}(x) = \dfrac{1}{2 \, \beta^4} \, \text{ln} \left( 1 + \exp(\beta^2 \, x)\right)^2 \, .
    \label{eq:SoftplusSquared}
\end{equation}
It corresponds to a squared version of \Cref{eq:Softplus} with a squared sharpness parameter for convexity and a divisor of two for a Softplus-like first derivative.
Note that its unbounded output and \ac{ReLU}-like rectification are prone to exploding outputs and gradients, especially in deep architectures.

The sharpness parameter $\beta$ can be included as a learnable parameter.
The preceding work has used a single parameter for the whole \ac{ANN}~\cite{Asad2022, Fleres2025}.
While some studies propose per-neuron parameterization~\cite{He2015,Apicella2021}, this work follows its reference and deploys the SoftplusSquared function in \Cref{eq:SoftplusSquared} with a single global parameter.

\subsubsection{Constraining parameters}

Another requirement of the \ac{ICNN} architecture are non-negative weights in all linear layers but the first.
Wrapping an appropriate function around the weights constitutes a possible implementation, e.g., absolute, square, exponential,  Softplus, or clipping.
In line with previous work~\cite{Asad2022, Fleres2025}, this study deploys the Softplus function for that purpose.
Hence, a constrained weight $\bar{w}_{ij}^{(l)}$ of layer $l$ is computed from its unconstrained value $w_{ij}^{(l)}$ as
\begin{equation}
    \bar{w}_{ij}^{(l)} = \text{Softplus}(w,\alpha) = \dfrac{1}{\alpha^2} \, \text{ln} \left( 1+\exp(\alpha^2 \, w_{ij}^{(l)})\right) \, .
\end{equation}
The squared sharpness parameter $\alpha$ ensures positive numbers and is included as a single learnable parameter~\cite{Asad2022,Fleres2025}.

\subsubsection{Correction terms for consistency}

The original work on \ac{ANN}-based constitutive modeling introduced a correction term to ensure consistent materials, i.e., materials that map zero strain to zero stress~\cite{Asad2022}.
In case of hyperreduction~\cite{Fleres2025}, the correction term ensures zero restoring forces at zero displacement.

The idea is to not use the \ac{ANN}'s prediction directly, but to modify it by a linear correction term.
This correction term is linear for the energy and, thus, constant for the force.
As the force must be zero at zero input, the offset is chosen to be the \ac{ANN}'s predicted force at zero input.
As a consequence, the predicted strain energy only deviates by a constant.
The resulting function of \ac{ICNN} and correction term is referred to as \ac{PANN} because it includes physical knowledge for zero input by design.

This work proposes an additional constant energy offset to ensure consistent energy, i.e., zero strain energy at zero displacement.
Similarly to the force correction term, the energy correction subtracts the prediction for a zero input.
As this work proposes second-order Sobolev training and thus, also learning the strain energy, the correction ensures physical behavior and removes the final layer's bias terms.
The constant offset does not interfere with the force correction.

In summary, the predictions of the \ac{PANN} for energy $\hat{e}$, force $\hat{\vi{f}}_r$, and stiffness $\hat{\vi{K}}_r$ based on the underlying \ac{ICNN} read
\begin{alignat}{3}
    & \hat{e}(\vi{x}_r) && = \hphantom{\partialfrac{}{\vi{x}_r}} \text{PANN}(\vi{x}_r) \\
    & &&= \hphantom{\partialfrac{}{\vi{x}_r}} \text{ICNN}(\vi{x}_r) - \partialfrac{}{\vi{x}_r} \, \text{ICNN} \, (\vi{0}) \cdot \vi{x}_r - \text{ICNN} \, (\vi{0}) \notag \\ 
    & \hat{\vi{f}_r}(\vi{x}_r) && = \partialfrac{}{\vi{x}_r} \, \text{PANN} \, (\vi{x}_r) \\
    & &&= \partialfrac{}{\vi{x}_r} \, \text{ICNN} \, (\vi{x}_r) - \partialfrac{}{\vi{x}_r} \, \text{ICNN} \, (\vi{0}) \notag \\ 
    & \hat{\vi{K}_r}(\vi{x}_r) && = \ppartialfrac{}{\vi{x}_r} \, \text{PANN} \, (\vi{x}_r) \\
    & &&= \ppartialfrac{}{\vi{x}_r} \, \text{ICNN} \, (\vi{x}_r) \, . \notag
\end{alignat}
Note that these correction terms must be also used in training and require a modified forward pass for prediction.

\subsubsection{Physical initialization of quadratic pass-through}

Another feature proposed in the reference work~\cite{Asad2022,Fleres2025} is a pass-through layer that directly connects the input to the output as shown in \Cref{fig:ICNN_architecture}.
The specific function has been chosen quadratic and reads
\begin{equation}
    f(\vi{x}_r) = \vi{x}_r^\top \, \vi{A} \, \vi{x}_r \, .
    \label{eq:quadratic_passthrough}
\end{equation}
To ensure convexity, $\vi{A}$ must be positive semi-definite. 
To guarantee this property by design, $\vi{A}$ can be constructed an invertible matrix as $\vi{A}=\Tilde{\vi{A}}^\top \, \Tilde{\vi{A}}$~\cite{Asad2022}.
Its coefficients are initialized randomly and constitute trainable parameters.

This work proposes an alternative initialization that incorporates a priori knowledge.
A prerequisite is the observation that \Cref{eq:quadratic_passthrough} exhibits the same structure as the expression for strain energy of a linear model where $\vi{A}$ is half of the stiffness matrix.
Therefore, we propose to initialize $\vi{A}$ as half of the tangent stiffness matrix for zero displacement, $\mi{K}_{r,0}$, potentially accelerating training.
To incorporate this initialization into the layer's definiteness-ensuring construction as $\vi{A}=\Tilde{\vi{A}}^\top \, \Tilde{\vi{A}}$, we utilize its eigendecomposition
\begin{equation}
    \mi{K}_{r,0} = \mi{Q} \, \mi{\Lambda} \, \mi{Q}^\top \, ,
\end{equation}
where $\mi{K}_{r,0}$ is the reduced tangent stiffness matrix at zero displacement, $\mi{Q}$ an orthogonal matrix containing its eigenvectors, and $\mi{\Lambda}$ a diagonal matrix of corresponding eigenvalues.
By choosing
\begin{equation}
    \Tilde{\vi{A}} = \frac{1}{\sqrt{2}} \, \left( \mi{Q} \, \mi{\Lambda}^{\circ \frac{1}{2}} \right)^\top \, ,
\end{equation}
where $^{\circ \frac{1}{2}}$ denotes the elementwise root or Hadamard root, we ensure that
\begin{equation}
    \begin{aligned}
        \Tilde{\vi{A}}^\top \, \Tilde{\vi{A}} 
        &= \frac{1}{2} \, \left(\mi{Q} \, \mi{\Lambda}^{\circ \frac{1}{2}}\right) \left(\mi{Q} \, \mi{\Lambda}^{\circ \frac{1}{2}}\right)^\top \\
        &= \frac{1}{2} \, \mi{Q} \, \mi{\Lambda}^{\circ \frac{1}{2}} \, \left(\mi{\Lambda}^{\circ \frac{1}{2}}\right)^\top \, \mi{Q}^\top \\
        &= \frac{1}{2} \, \mi{Q} \, \mi{\Lambda} \, \mi{Q}^\top \\
        &= \frac{1}{2} \, \mi{K}_{r,0}
    \end{aligned}
\end{equation}
and thus
\begin{equation}
\begin{aligned}
    f(\vi{x}_r) =& \vi{x}_r^\top \, \vi{A} \, \vi{x}_r \\
            =& \vi{x}_r^\top \, \Tilde{\vi{A}}^\top \, \Tilde{\vi{A}} \, \vi{x}_r \\
            =& \frac{1}{2} \, \vi{x}_r^\top \, \mi{K}_{r,0} \, \vi{x}_r \, ,
\end{aligned}
\end{equation}
which corresponds to the strain energy of a linear system.

\subsubsection{Data preprocessing}
\label{sssec:data_preprocess}

Another aspect is to scale data as it often leads to better performance~\cite{Shanker1996,LeCun1998}.
The reference work~\cite{Fleres2025} stated that min-max normalization cannot be deployed as it compromises convexity.
However, common scaling techniques such as normalization or standardization do not affect convexity of an \ac{ICNN}: 
These transformations can be interpreted as linear functions, similarly to linear layers with weights and biases.
For both standardization and normalization, all weights are non-negative.
Furthermore, composing one of these these linear transformations for the inputs with the first regular linear layer forms a single linear layer.
Therefore, it does not harm convexity.

However, another difficulty arises due to Sobolev training:
The different types of outputs cannot be scaled independently as they are related as derivatives.
Hence, if the \ac{ANN} was to learn scaled outputs, the scaling must be considered for the corresponding derivatives using the chain rule~\cite{Kissel2020}.
There are heuristics~\cite{Kissel2020} how to choose appropriate scaling, often aiming for similar orders of magnitude, but this branch is still a subject of research.

This work investigates the effect of standardized inputs to potentially enhance the \ac{ANN}'s performance, especially for accuracy and generalization.
Outputs remain unscaled as their scaling is ongoing research.
In detail, the input is standardized per component, as the reduced states for a \ac{POD} basis often feature different scales.
Standardization is included as a frozen linear layer placed before the original \ac{ANN}.
This choice ensures convenient handling as no outputs or derivatives must be transformed.
Furthermore, it also simplifies deployment as the \ac{ANN} forms a stand-alone solution.

Another aspect is the data's division into sets for training, validation, and testing.
These are used for adjusting the \ac{ANN}'s weights, to indicate performance during training and to identify an initialization's best model, and to assess its performance on completely unseen data, respectively. %
Note that validation data affect training as it determines which model performs best during training and is therefore saved.
Note that train and validation data are randomly sampled for each initialization while the test set is the same for all models.

\subsubsection{Loss function}
\label{sssec:loss}

The loss function assesses an \ac{ANN}'s performance and should be minimized during training.
A well-established choice for regression problems is the \ac{MSE} error~\cite{Asad2022}, or the \ac{hMSE}~\cite{Fleres2025}. 
Following the reference work, we deploy the latter, which in case of reduced forces reads
\begin{equation}
    \mathcal{L}_F = \dfrac{1}{2} \, \dfrac{1}{N_S} \,\dfrac{1}{r} \, \sum_{i=1}^{N_S}\, \sum_{k=1}^{r} \left( \bigl( \vi{\hat{f}}_{r,i} \bigr)_k- \bigl( \vi{f}_{r,i} \bigr)_k \right)^2 \, ,
    \label{eq:loss_f}
\end{equation}
where $r$ is the reduced dimension and $N_S$ the number of samples, i.e., $i$ selects the sample and $k$ the component.
The sum of squared component-wise differences can also be expressed as the vector's squared $L_2$-norm.

If the target is a vector and its components vary significantly in magnitude, component-wise scaling by the respective standard deviation can improve predictions of small components~\cite{Avrutskiy2020,Asad2022}. 
However, note that typical sampling of \ac{FEM} models, i.e., typical load cases, and reduction via \ac{POD} might produce components close to zero.
In this case, dividing by their negligible element-wise standard deviation significantly magnifies losses by several orders of magnitude and thus, disturbs training.
Therefore, this work does not utilize component-wise scaling by the respective standard deviation.

Another major aspect of this work is the extension from single to multiple losses due to second-order Sobolev training, i.e., training the function and its first two derivatives.
In contrast, the reference work~\cite{Fleres2025} has focused exclusively on the first derivative, which corresponds to the reduced force and constitutes the quantity of interest in hyperreduction.
However, the function itself and its Hessian, i.e., the strain energy and tangent stiffness matrix, are also available from commercial software and contain additional potentially powerful information.
To include all three quantities in the loss function, their respective \acp{hMSE} are computed analogously to \Cref{eq:loss_f}, i.e., 
\begin{align}
    \mathcal{L}_E &= \dfrac{1}{2} \, \dfrac{1}{N_S} \, \sum_{i=1}^{N_S} \left( \hat{e}_i - e_i \right)^2 \\
    \mathcal{L}_K &= \dfrac{1}{2} \, \dfrac{1}{N_S} \,\dfrac{1}{r^2} \, \sum_{i=1}^{N_S}\, \sum_{j=1}^{r} \, \sum_{k=1}^{r} \left( \bigl( \vi{\hat{K}}_{r,i} \bigr)_{jk} - \bigl( \vi{K}_{r,i} \bigr)_{jk} \right)^2 \, .
\end{align}

In order to incorporate all objectives in parameter updates, the three losses must be combined.
However, they might be of different orders of magnitude, of varying difficulty to learn, and might result in conflicting gradients for parameter updates.
Therefore, combining multiple losses poses a challenge on its own with a notable body of research.
Potential approaches include weighted sums~\cite{Pukrittayakamee2011,Avrutskiy2020,Feng2022c,Kim2022}, scheduled losses~\cite{Bouhlel2020,Kissel2020,Krishnapriyan2021}, dynamically weighted sums~\cite{Wang2023b}, and also aggregating and unifying the individual losses' parameter updates~\cite{Yu2020,Quinton2024}.
The following paragraphs describe several approaches that are later empirically studied within the numerical example.

The simplest approach to assemble a single loss function $\mathcal{L}$ is the sum of the individual losses, i.e.,
\begin{equation}
    \mathcal{L} = \mathcal{L}_E + \mathcal{L}_F + \mathcal{L}_K \, .
\end{equation}
Although this approach establishes a baseline for comparison, it does not take into account the challenges mentioned above and is therefore likely to produce poor results.

A straightforward modification to mitigate the different orders of magnitude is to deploy a weighted sum.
The weights can be chosen based on intuition~\cite{Feng2022c}, empirical studies~\cite{Kim2022}, standard deviations~\cite{Avrutskiy2020}, or the squared ratio of maximum values~\cite{Pukrittayakamee2011}.
A related idea is to train only the Hessian: if two functions have the same Hessian, the gradient can only differ by a constant, and the function only by linear and constant parts.
However, the correction terms for consistent energy and force eliminate these deviations.

A complementary concept is to schedule the individual losses, gradually including~\cite{Bouhlel2020} or excluding~\cite{Kissel2020} information.
The former aims for faster convergence and more accurate predictions of both function and gradient.
A similar approach has been proposed for \acp{PINN}~\cite{Krishnapriyan2021}, where the derivative-based loss is progressively included in consecutive trainings, significantly improving accuracy.
Therefore, the case study investigates a gradually increasing schedule for the Hessian.

Another potential problem are conflicting parameter updates of the individual losses.
This problem also occurs in multi-task learning where it is approached by aggregating the individual gradients.
Aggregation denotes merging the vectors via additional computations to minimize conflicts.
Notable examples of this technique include \ac{PCGrad}~\cite{Yu2020} and \ac{JD}~\cite{Quinton2024}.
The latter is the more recent algorithm and ensures a non-conflicting parameter update as a linear combination of the individual gradients. Furthermore, each gradient contributes proportionally to its norm.
This functionality is available in TorchJD~\cite{Quinton2024}, a python library for PyTorch.

However, the techniques described so far are not necessarily appropriate for objectives of different difficulties.
In other words: some losses may converge fast, so that more inert losses dominate the total loss, limiting improvements for the fast converging ones.
A possible remedy is the use of dynamically weighted losses~\cite{Wang2023b,Bischof2025} to maintain comparable orders of magnitude among the loss terms.
An example in the context of \acp{PINN} is loss balancing~\cite{Wang2023b}, which guarantees gradients of equal norm for the individual losses.
In detail, loss balancing scales each gradient to have the same norm as the sum of the individual gradients' norms.
The corresponding dynamic scaling factor $\hat{\lambda}$ for loss $i$ is given by
\begin{equation}
    \hat{\lambda}_i = \dfrac{\sum_j \norm{\nabla_\theta \mathcal{L}_j}}{\norm{\nabla_\theta \mathcal{L}_i}} \, ,
    \label{eq:loss_balancing_weight}
\end{equation}
where $\nabla_\theta \mathcal{L}_j$ denotes the gradient of loss $j$ with respect to the \ac{ANN}'s trainable parameters $\theta$.
A moving average with with $\alpha = 0.9$ smoothens the weights
\begin{equation}
    \lambda_{new} = \alpha \, \lambda_{old} + (1 - \alpha) \, \hat{\lambda}_{new} \, .
    \label{eq:MultiLoss_annealing_moving_avg}
\end{equation}

Note that the procedure amplifies low losses to the magnitude of dominant ones.
In case of an inert loss, a potential risk are parameter updates of too large magnitudes for the fast converging ones, corresponding to a learning rate chosen too high.
However, the concept can be easily adapted, e.g., to decrease high losses to the level of lower ones.
Therefore, this work proposes an alternative dynamic loss balancing scheme that aims for same scale losses.
The idea is to weight each loss based on the ratio of the average of all other losses to the total sum of losses
\begin{equation}
    \lambda_i 
    = \dfrac{\sum_{j, \, j \neq i} \mathcal{L}_j}{N_{L} - 1 } \, \dfrac{1}{ \sum_j \mathcal{L}_j } 
    = \dfrac{1}{N_{L} - 1 } \, \left( 1 - \dfrac{\mathcal{L}_i}{ \sum_j \mathcal{L}_j }\right) \, ,
\end{equation}
where $N_{L}$ is the number of loss terms.
With this balancing scheme, same scale losses lead to equal scaling factors.
In case of a dominating loss, it is scaled down to the level of the remaining ones.
Moreover, the sum of all scaling factors $\lambda_i$ equals one.
The numerical case study in \Cref{ssec:Results} investigates both the original loss balancing and the modified implementation.

\section{Numerical case study: hyperelastic beam}
\label{ssec:Results}

This section applies the two hyperreduction methodologies to a numerical case study.
A static nonlinear \ac{FEM} model serves as a reference and data source.
The specific model is reproduced from the reference work~\cite{Fleres2025} for comparability with a detailed description in \Cref{ssec:FEM_model}.
This model and its data are reduced via the \ac{POD} with results presented in \Cref{ssec:MOR_results}.
The reduced data are the basis for hyperreduction, either via \ac{TPWL} or a \ac{PANN}.
Before comparing these two methods, different aspects of the \ac{PANN} approach are empirically studied in \Cref{ssec:PANN_results}.
These studies comprise two parts.
The first part optimizes hyperparameters for a \ac{PANN} trained on force only, in line with reference work~\cite{Fleres2025}.
On that basis, the second part studies second-order Sobolev training with special emphasis on combining the individual losses.
After the best settings have been identified, \Cref{ssec:Hyperreduction_results} assesses the accuracy of the best performing \ac{PANN} model, a \ac{TPWL} model, and the \ac{PANN} of the reference work~\cite{Fleres2025}.

All \ac{FEM}-related work utilizes Ansys\textsuperscript{\textregistered} Academic Research Mechanical APDL, Release 2022 R2, and MKS units. 
The model is scripted in APDL and is publicly available on GitLab at \url{https://gitlab.gwdg.de/jade-hochschule/fms/2025_cmame}.
All subsequent steps are implemented in python, especially relying on PyTorch.
The corresponding scripts are available on GitLab, too.

\subsection{Finite element model}
\label{ssec:FEM_model}

\begin{figure}[!htb]
    \centering
    \includegraphics[width=0.5\linewidth]{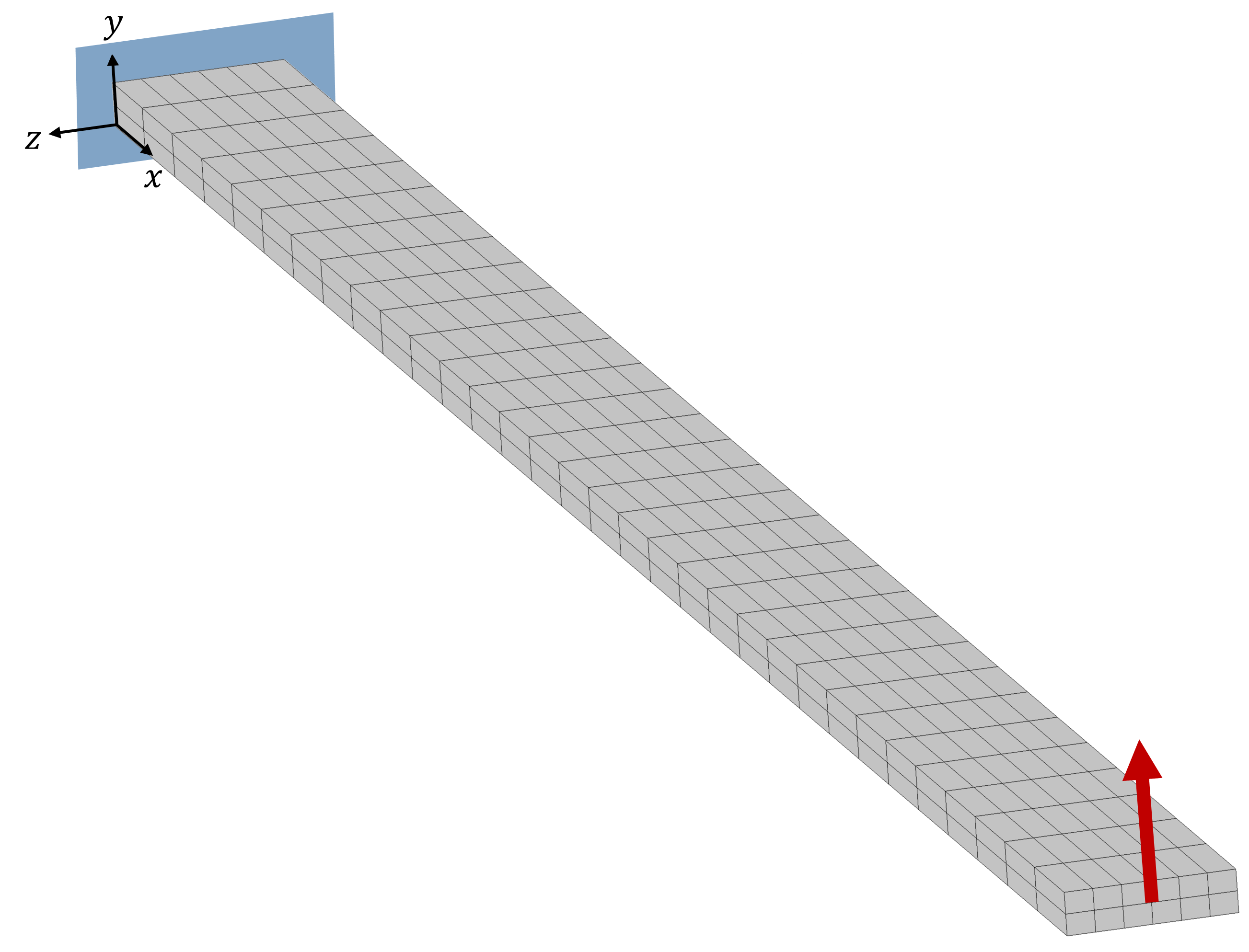}
    \caption{\ac{FEM} model of the numerical case study. A cantilever beam is fixed at one end and subject to a tip load at the free end.}
    \label{fig:FEM_beam_mesh}
\end{figure}

The \ac{FEM} model described here constitutes the reference and basis for subsequent \ac{MOR} and hyperreduction.
As this study builds on top of a reference work~\cite{Fleres2025}, its case study is reproduced for comparability.
This case study comprises a nonlinear static structural analysis of a cantilever beam subjected to an increasing tip load.
The cantilever beam is of dimension $160 \times 2 \times \qty{8}{\milli\meter}$ and consists of nearly-incompressible isotropic neo-Hookean material.
Its strain-energy potential is given by~\cite[ch.~7]{Rust2020}
\begin{equation}
    W = \dfrac{G}{2} \, \left( \bar{I}_1 - 3 \right) + \dfrac{K}{2} \left( J - 1 \right) \, ,
\end{equation}
where $G$ denotes the initial shear modulus, $K$ the initial bulk modulus, $J$ the volume ratio, and $\bar{I}_1$ a modified strain invariant.
$K$ can be expressed in terms of the Poisson's ratio $\nu$ and the initial Young's modulus $E$ as
\begin{equation}
    K = \dfrac{E}{3 \, (1 - 2 \nu)} \, .
\end{equation}
Similarly, $G$ can be computed as 
\begin{equation}
    G = \dfrac{E}{3} \, .
\end{equation}
This work uses $E=\qty{111}{\giga\pascal}$ and $\nu = 0.33$, which are chosen from the parameter ranges of the reference work~\cite{Fleres2025}.

\begin{figure}[!htb]
    \centering
    \includegraphics[width=0.75\linewidth]{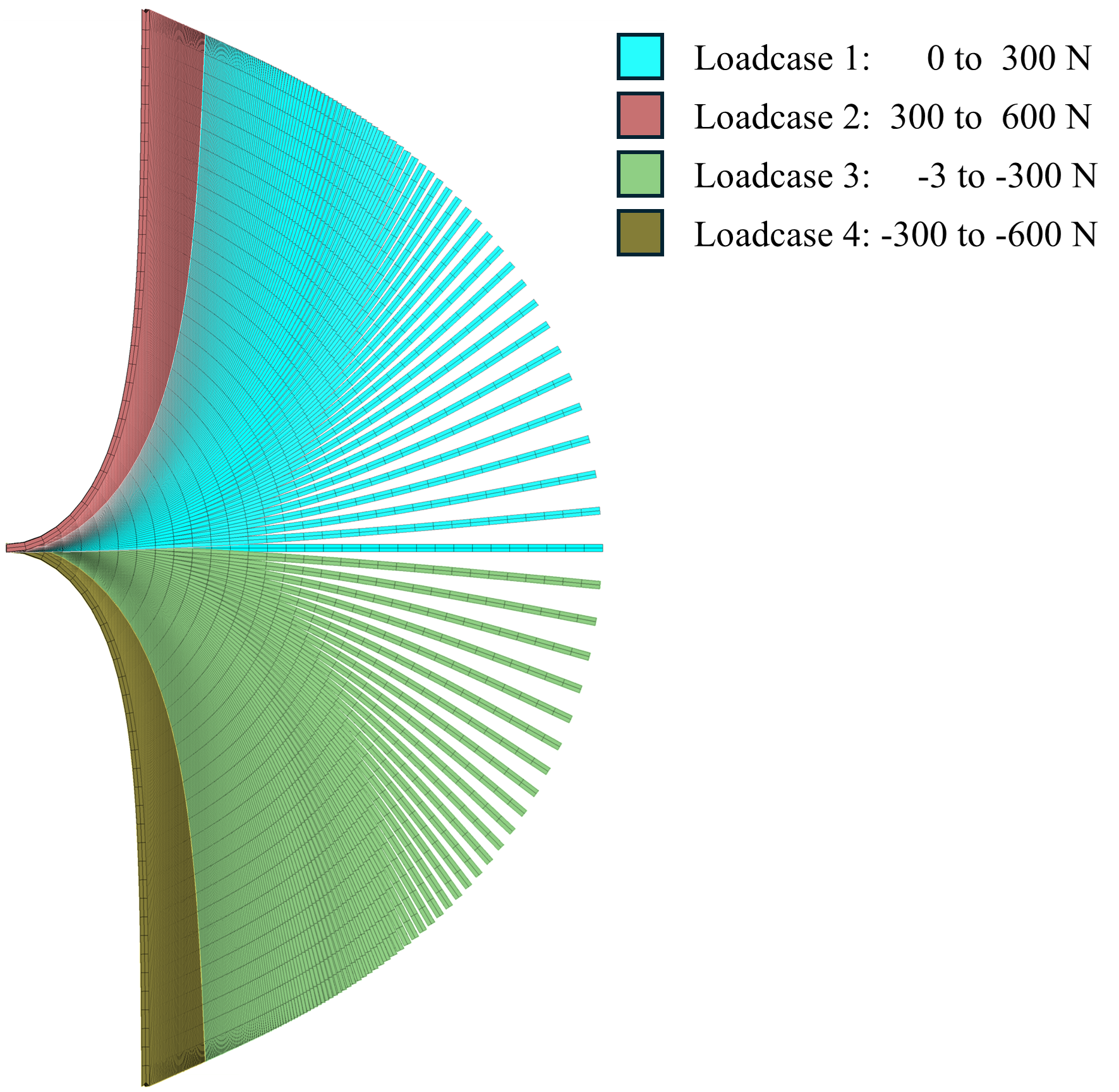}
    \caption{Side view of beam deformation for all load steps where color indicates the four load cases. Note that only the first load case colored in cyan is used for constructing or training the hyperreduced model while the other three constitute a challenging test.}
    \label{fig:FEM_beam_loadcases}
\end{figure}

The computational domain is spatially discretized into $32 \times 2 \times 6$ hexaedral elements with linear shape functions. 
\Cref{fig:FEM_beam_mesh} illustrates the computational mesh.
The discretization slightly differs from the reference, which uses $32 \times 2 \times 5$ elements.
This difference produces a node centered at the cantilever beam's free end, which will be used to apply the load and to measure its displacement.
The resulting mesh comprises $384$ elements and $693$ nodes with three degrees of freedom each.
The fixed end affects $21$ nodes, resulting in a total of $n=2016$ unknowns.
As in the reference, the load case comprises a point force in positive y-direction with a maximum amplitude of \qty{300}{\newton}, which builds up gradually in $100$ homogeneous load steps.
In order to critically assess the generalization capability of the hyperreduction methods, three more load cases are added.
These include increasing the force magnitude up to \qty{600}{\newton} and the same loading in negative y-direction.
\Cref{fig:FEM_beam_loadcases} illustrates all four load cases with respective beam deflections.
Note that the three additional load cases are only used for testing, i.e.,  are not part of constructing or training the hyperreduced model.
This choice of test data poses a rather hard challenge that likely causes high test errors.

The four load cases and the resting position result in $401$ load steps.
For each load step $i$, the sampled data include the state vector of displacements $\vi{x}_i \mydim{n}$, the strain energy $e_i \in \mathbb{R}$, the internal restoring force vector $\vi{f}_i \mydim{n}$, and the tangent stiffness matrix $\vi{K}_i \mydims{n}{n}$.
These data serve as a basis for subsequent dimensionality reduction and hyperreduction.

However, there are some differences between our setup and the reference work~\cite{Fleres2025}.
Most importantly, material parameters are fixed, leading to some differences.
The reason is the focus on hyperreduction, Sobolev training, and comparison to alternative hyperreduction via \ac{TPWL}.
This decision slightly simplifies the \ac{ANN} architecture as the parameter path is discarded.
Furthermore, the training goal is slightly simpler, as no parametric effects must be learned.
However, the setup in the reference work uses two material parameters, which have been sampled in a uniform $10 \times 10$ grid to study hyperparameters. 
Each point in this grid corresponds to a \ac{FEM} analysis with the same load case from \qtyrange{3}{300}{\newton}, therefore containing $100$ samples.
The training set includes \qty{80}{\percent} of the grid points, each with their full load case.
Therefore, the training in the reference work has used a total of \num{8000} samples.
In contrast, this work relies on only a single load case with $101$ samples, of which \qty{50}{\percent} are used in training.
Arguably, reducing the number of training samples by a factor of more than $150$ increases the challenge more than the slightly simpler goal and architecture reduce it.

\subsection{Results model order reduction}
\label{ssec:MOR_results}

\ac{MOR} and dimensionality reduction of the data follow the workflow described in \Cref{ssec:MOR}.
Therefore, \ac{POD} extracts a reduced basis from state vector snapshots via a \ac{SVD}.
The data reduced this way constitute the basis for training or constructing hyperreduced models in following subsections.

The snapshots matrix comprises all $401$ load steps.
This choice leaks data from the test data into training to some extent.
However, it minimizes deviations due to an unsuitable data representation.
Hence, deviations can be attributed more clearly as fewer effects contribute.
As this work focuses on the hyperreduction aspect via \acp{ANN}, this choice allows for more precise statements.
Moreover, extracting a reduced basis from training data only would require repeating the dimensionality reduction each initialization.

The reduced dimension is set to $r=4$, which corresponds to \qty{99.82}{\percent} cumulative energy.
Moreover, it matches the setup for studying hyperparameters deployed in the reference work~\cite{Fleres2025}.
A consistent reduced dimension is important when combining \ac{POD} and mean-based loss functions.
\ac{POD} introduces an descending order, which also reflects in the magnitudes of the corresponding reduced states.
Hence, an additional component with small magnitude might decrease the \ac{hMSE} loss.
However, a deviation arises because of different data.
The reference work's reduced basis includes parametric effects, which are omitted in our approach, whereas ours encompasses extended load cases.
Therefore, comparisons should focus on reduced quantities as the reduced basis does not contribute to their deviations. 

\begin{figure}[!htb]
    \centering
    \includegraphics[width=1\linewidth]{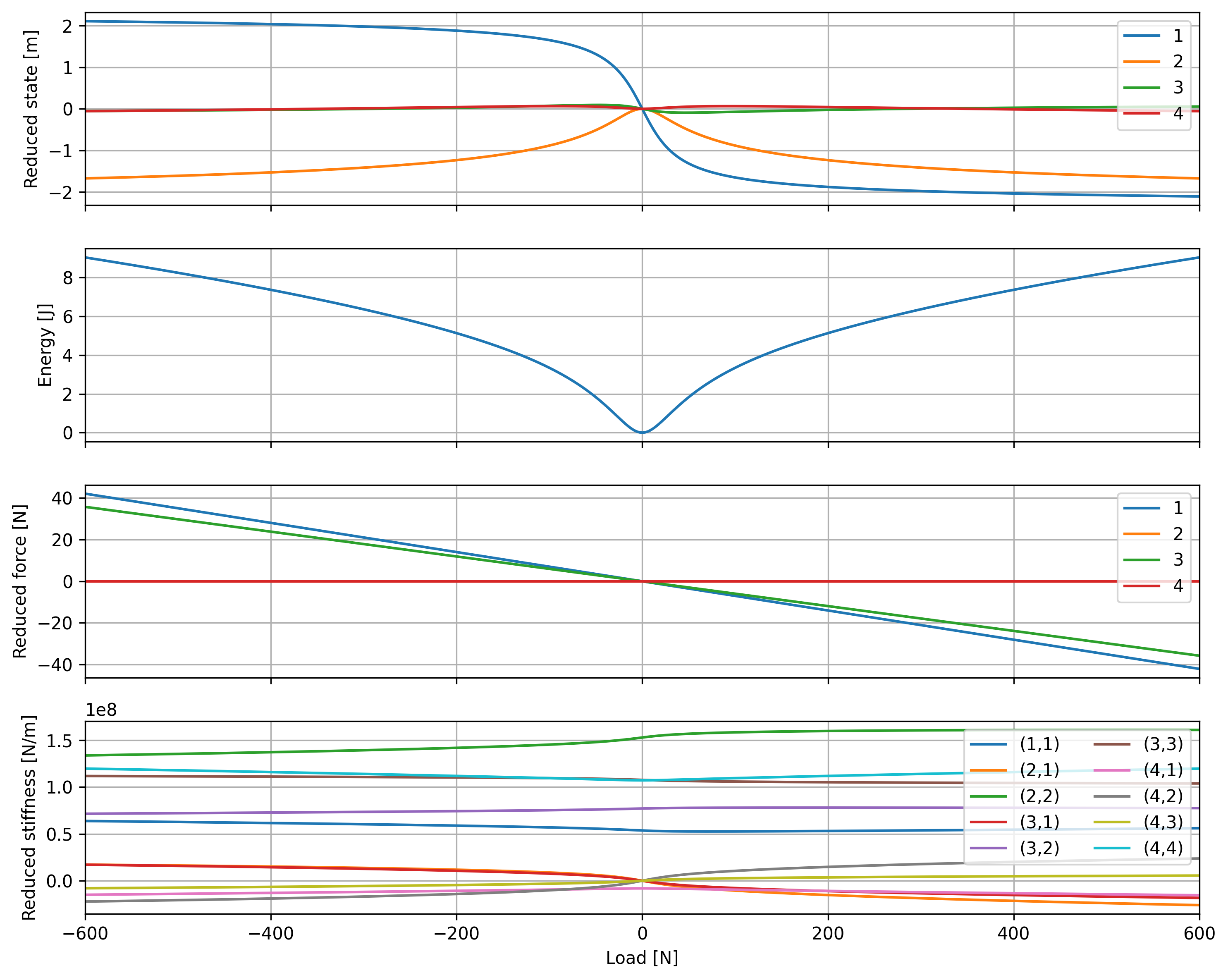}
    \caption{Reduced coordinates, strain energy, reduced force components, and reduced tangent stiffness components vs. load magnitude. Fewer matrix components are included due to symmetry. Note the stiffening behavior and the different orders of magnitude of energy, force, and tangent stiffness.}
    \label{fig:reduced_data}
\end{figure}

\Cref{fig:reduced_data} illustrates the reduced data, including the reduced state vector, the strain energy, the reduced internal force, and the reduced tangent stiffness matrix.
The reduced state shows a clear stiffening behavior.
After some initial deformation, the states quickly saturate.
For positive load direction, this effect is covered well by the training data.
Therefore, we expect low test losses for increasing the load from the reference load case by another \qty{300}{\newton}.
The stiffened regions potentially also help within an iterative nonlinear solution scheme because changes in external load lead to minor changes in displacement.
On the other hand, minor changes in the system's state can cause drastic changes in its internal force.
These two sensitivities must be considered when evaluating hyper-reduced models in \Cref{ssec:Hyperreduction_results}.

Note that the reduced coordinates vary greatly in magnitude.
This behavior results from \ac{POD} since it sorts the singular vectors by their contribution in explaining the data's variance in descending order.
As a result, the first reduced coordinate is the most relevant one, with subsequent ones progressively contributing less.
In order to promote learning all components, machine learning algorithms might benefit from preprocessing inputs to have the same scale.

Another observation is the difference in the magnitudes of energy, force, and stiffness.
The different scales are relevant in Sobolev training when combining the respective losses.
The reduced stiffness is several orders of magnitude higher.
In combination with a \ac{MSE}-based loss that penalizes extremes, this objective likely dominates the other two losses.

\subsection{Results physics-augmented neural network}
\label{ssec:PANN_results}

This subsection investigates hyperreduction based on \acp{PANN} via empirical studies in two parts:
The first part in \Cref{sssec:ANN_results_base} optimizes the \ac{PANN} architecture and hyperparameters.
While the reference work~\cite{Fleres2025} already studied some hyperparameters, it did not include the learning rate, which is arguably the most important hyperparameter~\cite{Bengio2012}.
Moreover, it considered only a single initialization per variant.
Therefore, the results do not reflect the inherently stochastic nature of \acp{ANN} with random parameter initializations and convergence to different local minima.
This work fills the gap and studies the learning rate with different architectures, initializing each variant $10$ times. 
In addition, the effects of standardizing input data and physically initializing the quadratic pass-through layer are investigated.
Other improvements include splitting data into three instead of two sets, effectively adding a validation set.
This set allows for monitoring performance during training and to checkpoint the best model.
A related change is to select the best hyperparameters based on validation loss, only using the test set for a final assessment and not for model selection.
Once optimal hyperparameters have been identified, the second part in \Cref{sssec:ANN_results_Sobolev} explores second-order Sobolev training with special emphasis on balancing multiple losses.
The best performing \ac{ANN} is then deployed for hyperreduction and compared with \ac{TPWL} and the reference work~\cite{Fleres2025}.

Every \ac{ANN} setup is trained with the following settings:
The \ac{ANN}'s parameters are adjusted by the \ac{Adam} optimizer~\cite{Kingma2014} due to its robustness and for consistency with the reference work~\cite{Fleres2025}.
For the latter reason, all weights are initialized following the uniform Glorot scheme~\cite{Glorot2010} and all bias terms as zeros.
Note that initializations are often tailored to the type of activation function~\cite{Glorot2010,He2015} or for a special architecture~\cite{Cranmer2020,Hoedt2023}.
To mitigate the effect of non-optimized initializations, each model is trained for generous \num{100000} epochs.
Initializing each variant $10$ times reduces the effect of randomness and ensure reliable results.
A fixed random seed enhances reproducibility and comparability between different variants.
For each initialization, the model that achieves the best validation loss during training is saved and deployed for testing.
As detailed in \Cref{sssec:data_preprocess}, the data is split into three sets: training, validation, and testing.
The $101$ samples that correspond to loading from \qtyrange{0}{300}{\newton} are randomly split per initialization \num{50}/\qty{50}{\percent} into training and validation data.
The training data form a single batch due to the low number of samples.
The test data are chosen to assess a model's capability for generalization and extrapolation, in contrast to reference work where test data only assess interpolation.
In detail, the test data comprises the remaining three loading cases, i.e., $300$ samples.
Especially loading in the opposite direction differs significantly and is not covered in the training data.
Therefore, large test losses are expected.

\subsubsection{Study I: architecture and hyperparameters}
\label{sssec:ANN_results_base}

The first study builds on the reference work~\cite{Fleres2025} and aims to find optimal hyperparameters and assess enhancements.
These enhancements include the layer for standardizing inputs and the physical initialization of the quadratic pass-through layer. 
Each variant is initialized $10$ times to account for the inherently stochastic nature.
Three architectures with varying depths but constant number of hidden neurons are studied.
Furthermore, three logarithmically spaced learning rates are investigated.
\Cref{tab:hyperparameters} describes the specific parameters and their values.
In addition, the physical initialization of the quadratic pass-through layer is studied, but only for each of the three architectures with a fixed learning rate of \num{1e-3} and unprocessed inputs.
This setup results in $21$ variants to be studied.
Each variant is trained only on its first derivative, i.e., the reduced force, in alignment with the reference work~\cite{Fleres2025}.
In contrast to the reference work, our study also deploys a validation set, enabling extended monitoring during training and checkpointing.
A checkpoint saves the model state with the best validation loss during training instead of the final state.
Another difference of this work is that the best model is selected based on validation loss, with test loss used solely as a final indicator of accuracy.

\begin{table}[h]
    \centering
    \begin{threeparttable}  
    \caption{Hyperparameters and their values.}
    \label{tab:hyperparameters}
    \begin{tabularx}{0.75\textwidth}{ll}
    \toprule 
    Aspect & Values \\
    \midrule 
    Architecture & $3 \times 50$, $2 \times 75$, $1 \times 150$\\
    (Number of parameters) & (5822), (6472), (922) \\
    Learning rate & \num{1e-2}, \num{1e-3}, \num{1e-4} \\
    Input scaling & None, standardized \\
    \bottomrule 
    \end{tabularx}
    \end{threeparttable}
\end{table}

\begin{figure}[!htb]
    \centering
    \includegraphics[width=1\linewidth]{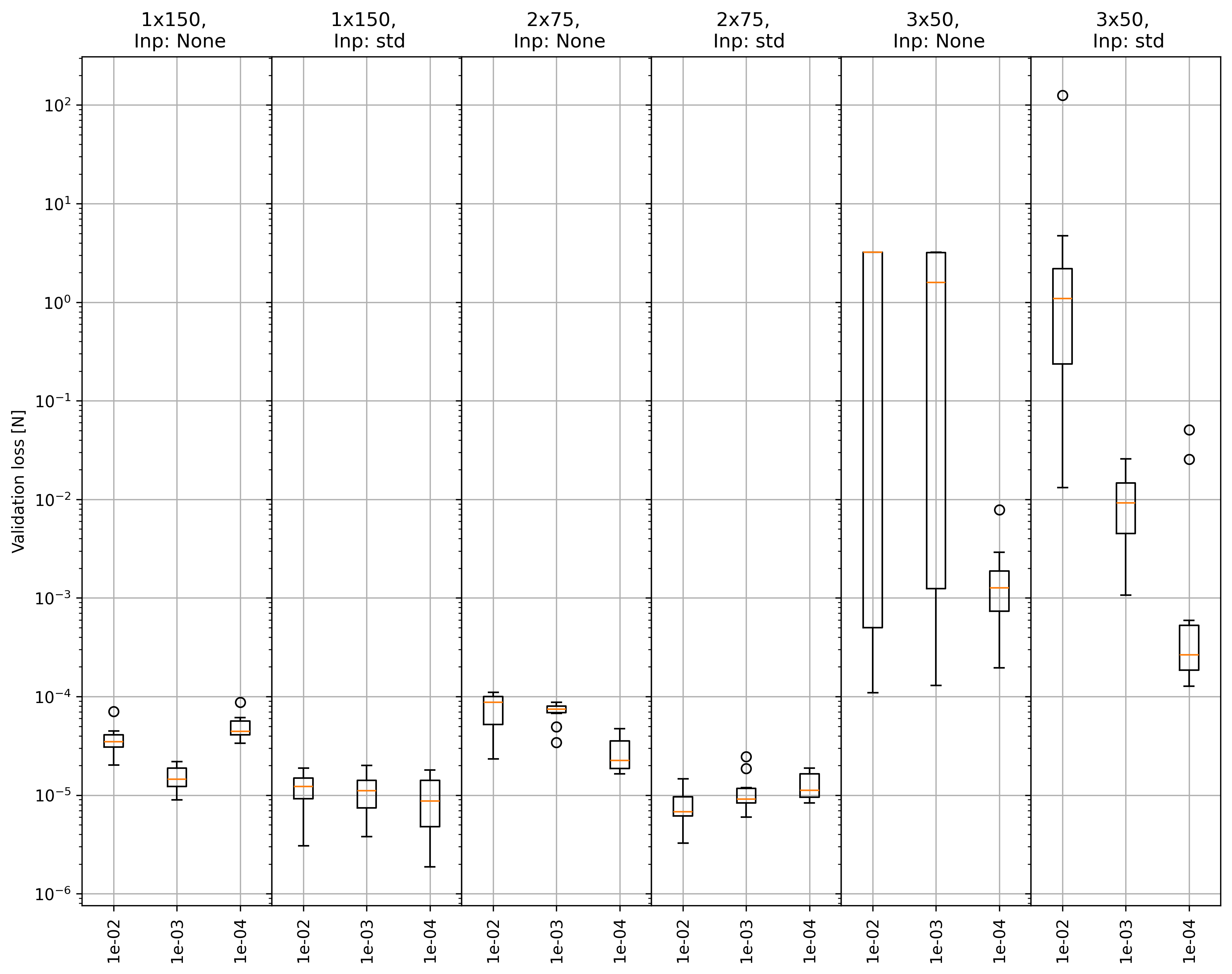}
    \caption{Statistical evaluation of model performances with respect to the validation error: each box summarizes the $10$ initializations per variant. Each of the six subplots corresponds to a combination of architecture and input scaling, while the three learning rates constitute the x-axis per subplot.}
    \label{fig:ANN_performance_validation}
\end{figure}

\Cref{fig:ANN_performance_validation} summarizes the different variants' best validation losses.
For each initialization of a specific variant, the model with the currently lowest validation loss is checkpointed.
Hence, its validation loss is the best value achieved during training.
One such loss constitutes a data point in the graph.
The $10$ validation losses per variant form a box plot.
There are $18$ boxes that correspond to all combinations of the hyperparameters in \Cref{tab:hyperparameters}.
The number of parameters per architecture is computed according to \Cref{app:N_par}.
The six subplots correspond to the three architectures with either unprocessed or component-wise standardized inputs.
Each subplot contains results for the three learning rates.
The following paragraphs compare to the reference work and extract concrete findings.

The reference work~\cite{Fleres2025} focuses on the $3 \times 50$-architecture with a learning rate of \num{1e-3} and unprocessed inputs.
This variant performed best in the respective hyperparameter sensitivity study, reaching a train loss of \qty{6e-5}{\newton} and a test loss of \qty{3e-4}{\newton}.
In our work, this variant's best initialization reaches a train loss \qty{1.34e-4}{\newton}, a validation loss of \qty{1.3e-4}{\newton}, and a test loss of \qty{1.56e2}{\newton}.
Therefore, the train loss is \qty{123}{\percent} higher in our reproduced result.
This deviation can be explained by slightly different setups, although consistency has been ensured wherever possible.
Most notably, this work discards parametric capabilities.
As elaborated in \Cref{ssec:FEM_model}, this difference slightly simplifies the \ac{ANN} architecture, but drastically reduces the amount of data.
In detail, the \ac{ANN} comprises \num{5.82e3} instead of \num{7.18e3} parameters as layers associated with material parameters become redundant.
Moreover, the number of samples for training decreases from \num{8000} to \num{51}.
Therefore, we expect the challenge to outweigh the simplification by a significant margin.
Hence, we attribute the increase in training loss to the more challenging setup of this work.
The other types of loss can hardly be compared because the test sets are fundamentally different.
Moreover, the reference work did not monitor performance during training via a validation set.
Further differences are multiple initializations, twice as many epochs, and checkpointing the model with the lowest validation loss instead of after the last epoch.
All these aspects potentially decrease losses, ensuring a fair or even optimistic representation of the reference work.

Considering all results shown in \Cref{fig:ANN_performance_validation}, the influence of the hyperparameters can be evaluated, and an optimal variant can be identified.
One global trend is standardized inputs improving the medians for all architectures.
For the $1\times150$ and $2\times75$-architectures, standardization also increases accuracy of a variant's best initialization.
Another aspect is the learning rate, to which most variants are rather robust.
An exception are the $3\times50$-architectures using standardized inputs, which clearly favor low values.
Variants with low learning rates might not have reached convergence, potentially causing underfitting and underestimating performance.
However, \num{100000} epochs are twice the value used in the reference work~\cite{Fleres2025} and are chosen to mitigate this effect.
The next observation is the good performance of shallow architectures in spite of a smaller capacity, indicating not yet optimal hyperparameters or full convergence for the deeper architectures.
Hence, tuning hyperparameters more extensively might reveal even more promising variants.
The lowest validation losses are achieved by variants with standardized inputs.
Promising candidates include the $1 \times 150$-architecture with all learning rates but slightly favoring lower ones, and the $2 \times 75$-architecture with a learning rate of \num{1e-2}.

\begin{figure}[!htb]
    \centering
    \includegraphics[width=1\linewidth]{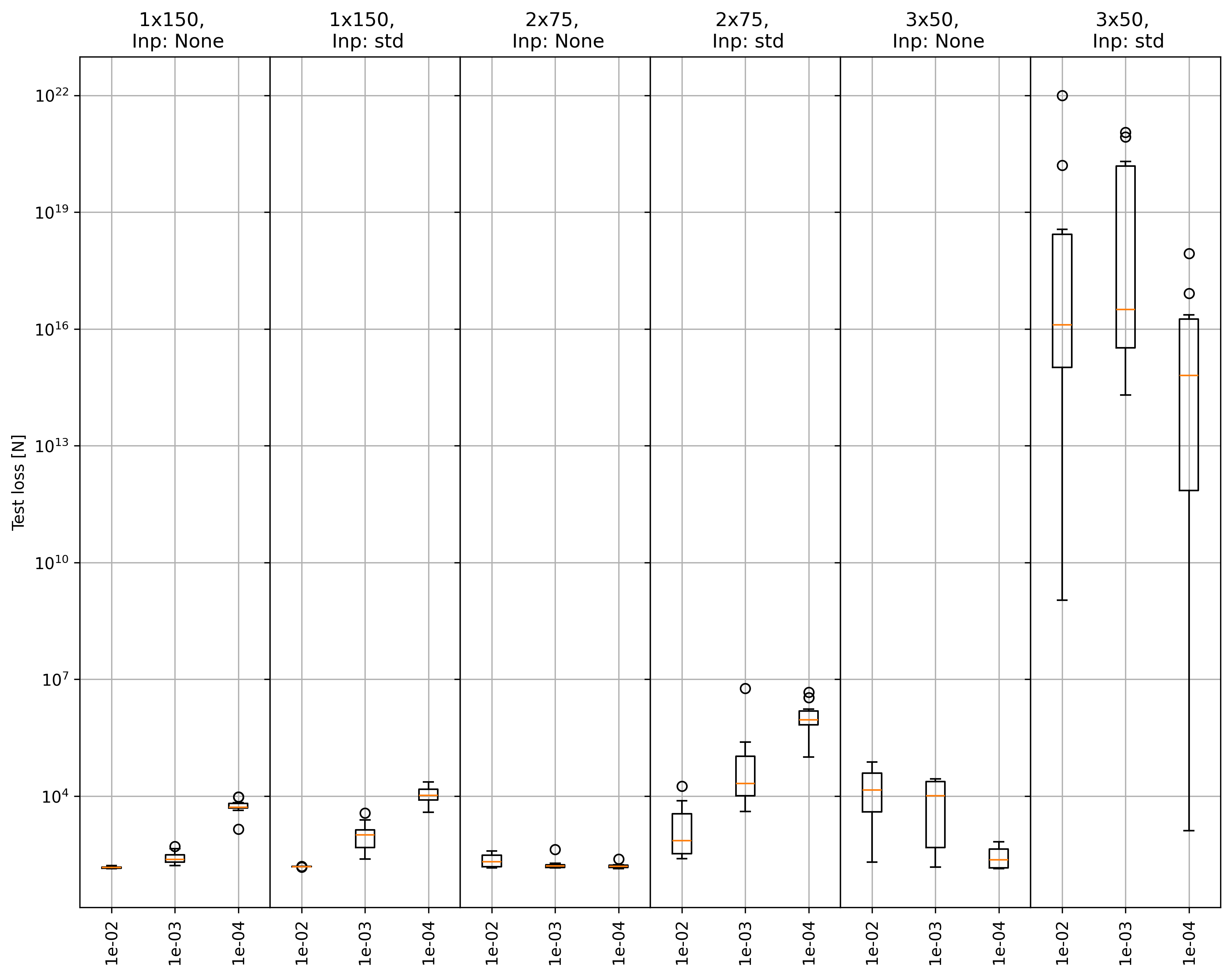}
    \caption{Analogous plot to \Cref{fig:ANN_performance_validation} to statistically summarize variants, but for the test error instead of validation error. Note that the test data exclusively comprises extrapolation load cases, posing a significant challenge.}
    \label{fig:ANN_performance_Val}
\end{figure}

\Cref{fig:ANN_performance_Val} resembles \Cref{fig:ANN_performance_validation} but contains test losses instead of validation losses.
Test data serve as a final evaluation, but must not be considered in model selection.
Note that the test data exclusively comprise extrapolation load cases.
These are not covered by train and validation data, leading to higher test losses.
The most prominent observation is the poor performance of variants with a $3 \times 50$-architecture and standardized inputs.
These add $15$ orders of magnitudes to the scale that are not required for any other variant.
Standardizing inputs also slightly worsens the performance of $2 \times 75$-architectures and negligible so for the $1 \times 150$-architectures.
The best performances are achieved by the two $1 \times 150$-architectures with a learning rate of \num{1e-2} and by $2 \times 75$-architectures with unprocessed inputs.
However, other variants reach similar accuracy, but only in outliers and far less consistently.

\begin{figure}[!htb]
    \centering
    \includegraphics[width=0.45\linewidth]{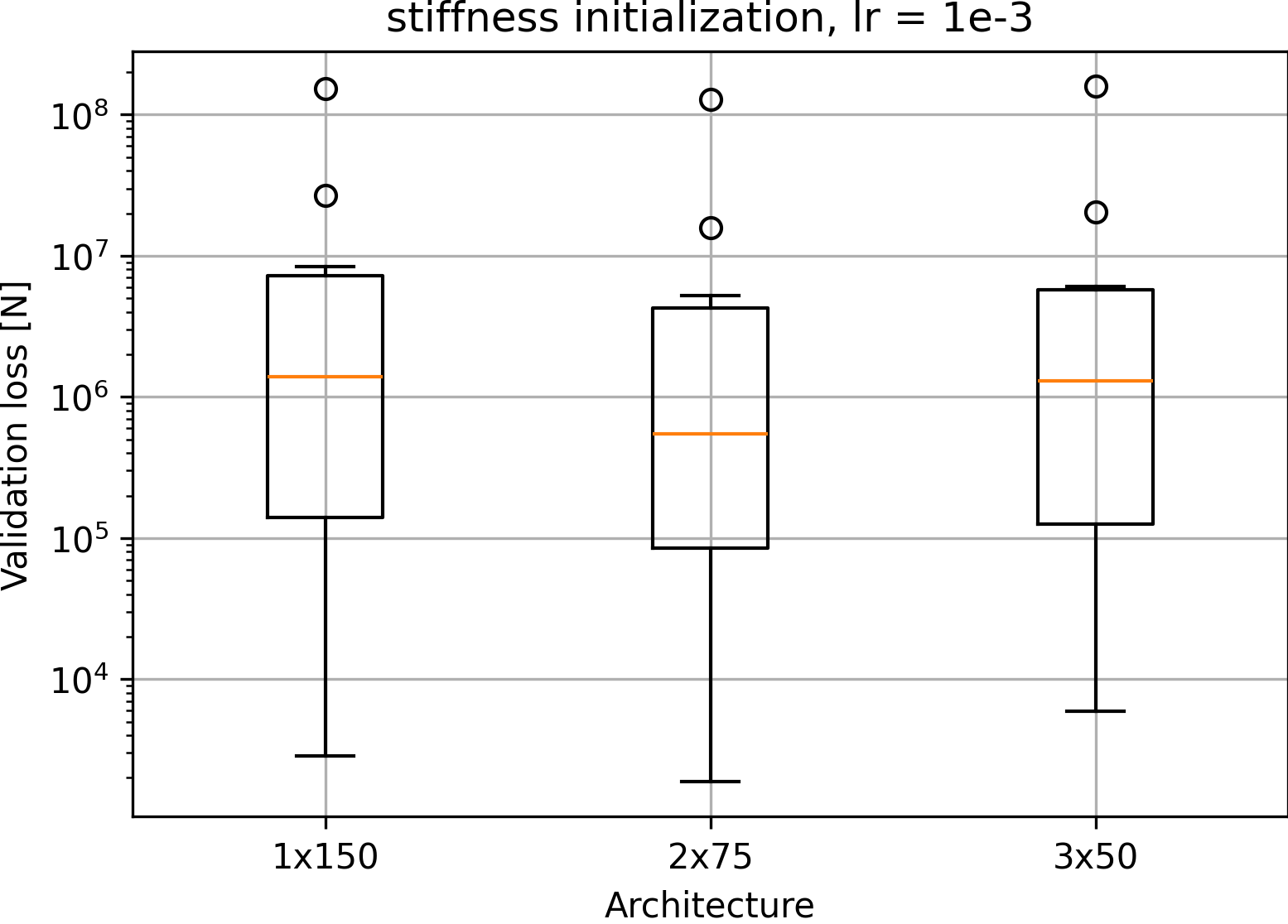}
    \hfill
    \includegraphics[width=0.45\linewidth]{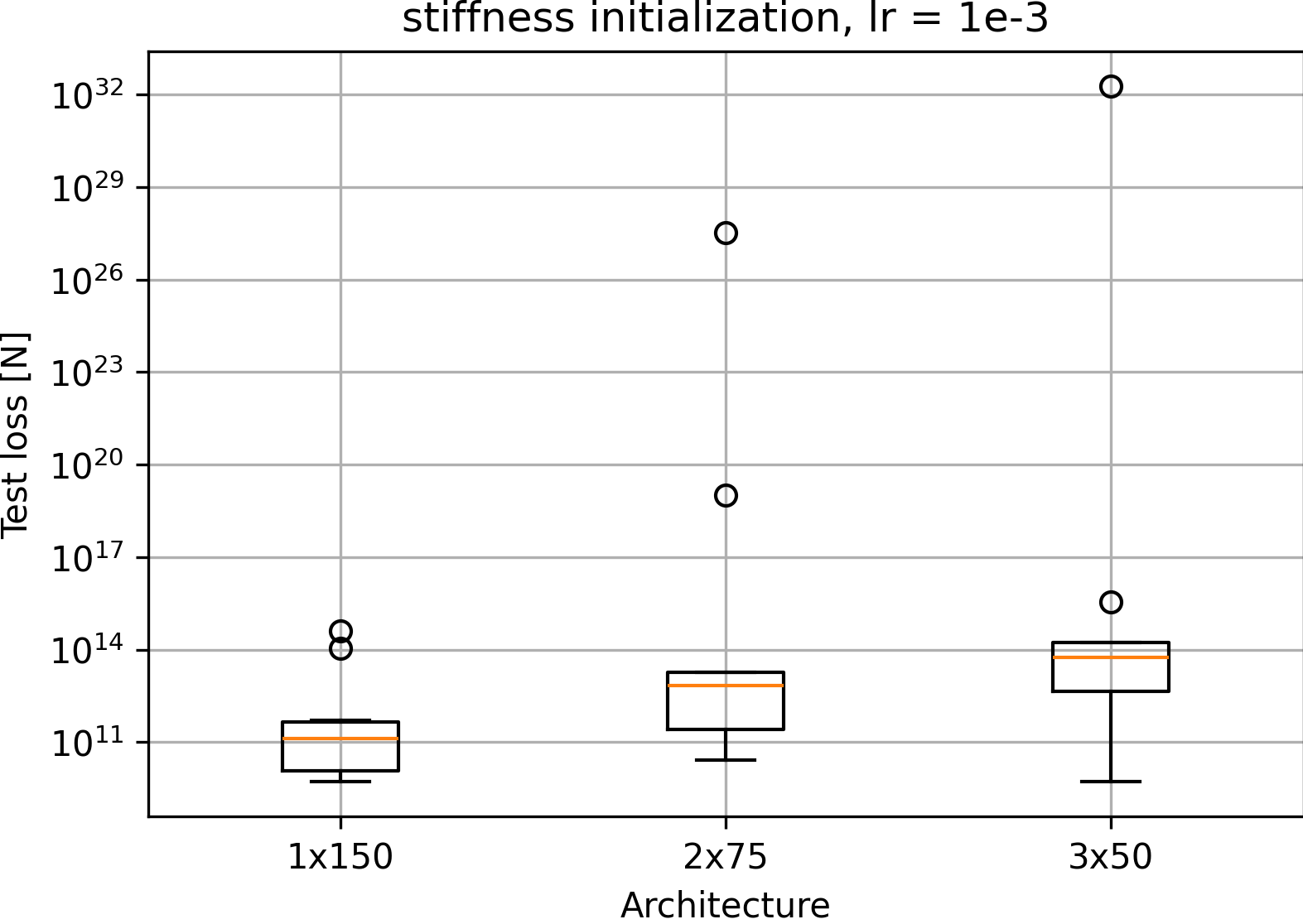}
    \caption{Statistical summary of validation and test losses of architectures where the quadratic pass-through layer has been initialized as the tangent stiffness matrix at rest. This option results in poor accuracy is thus, is discarded.}
    \label{fig:ANN_performance_stiffness}
\end{figure}

Another aspect to be studied is to initialize the quadratic layer's weight matrix as the tangent stiffness matrix at rest.
Therefore, this pass-through layer resembles computing the strain energy of a linear model.
\Cref{fig:ANN_performance_stiffness} presents the statistical summaries of validation and test losses for the corresponding variants.
All losses are several orders of magnitude worse than without this feature.
As a result, this concept of including a priori knowledge is discarded for future investigations.

The overall goal of this first study is to find a specific architecture and a set of hyperparameters as a basis for further studies with Sobolev training.
While $1 \times 150$-architectures seem promising regarding their losses and lower number of parameters, their second derivatives have only minimal representation power.
The reason is that a single layer of this activation function has a predominantly constant second derivative.
Therefore, one hidden layer does not suffice for learning the tangent stiffness matrix.
A promising alternative variant is found in the $2 \times 75$-architecture with a learning rate of \num{1e-2} and standardized inputs.
Its performs well for validation, achieving the lowest median from all variants, and also for test loss.
Its best initialization achieves a training loss of \qty{2.35e-06}{\newton}, a validation loss of \qty{3.27e-06}{\newton}, which is similar to the overall best validation loss of \qty{1.88e-06}{\newton}, and a test loss of \qty{2.47e2}{\newton}.
The train loss poses the most comparable quantity and surpasses the reference work with \qty{6e-05}{\newton} by a factor of $25$.
For that reason, this variant is chosen as the basis for Sobolev training.

In contrast, the best-performing model from the hyperparameter sensitivity study of the reference work achieves a training loss of \qty{3e-05}{\newton}.
Hence, the efforts of this work increased the training accuracy by one order of magnitude.
Note that there are additional differences, most dominantly due to discarded parametric capabilities.
However, as elaborated in \Cref{ssec:FEM_model}, these differences are expected to worsen accuracy for models in our setup.

\begin{figure}[!htb]
    \centering
    \includegraphics[width=0.75\linewidth]{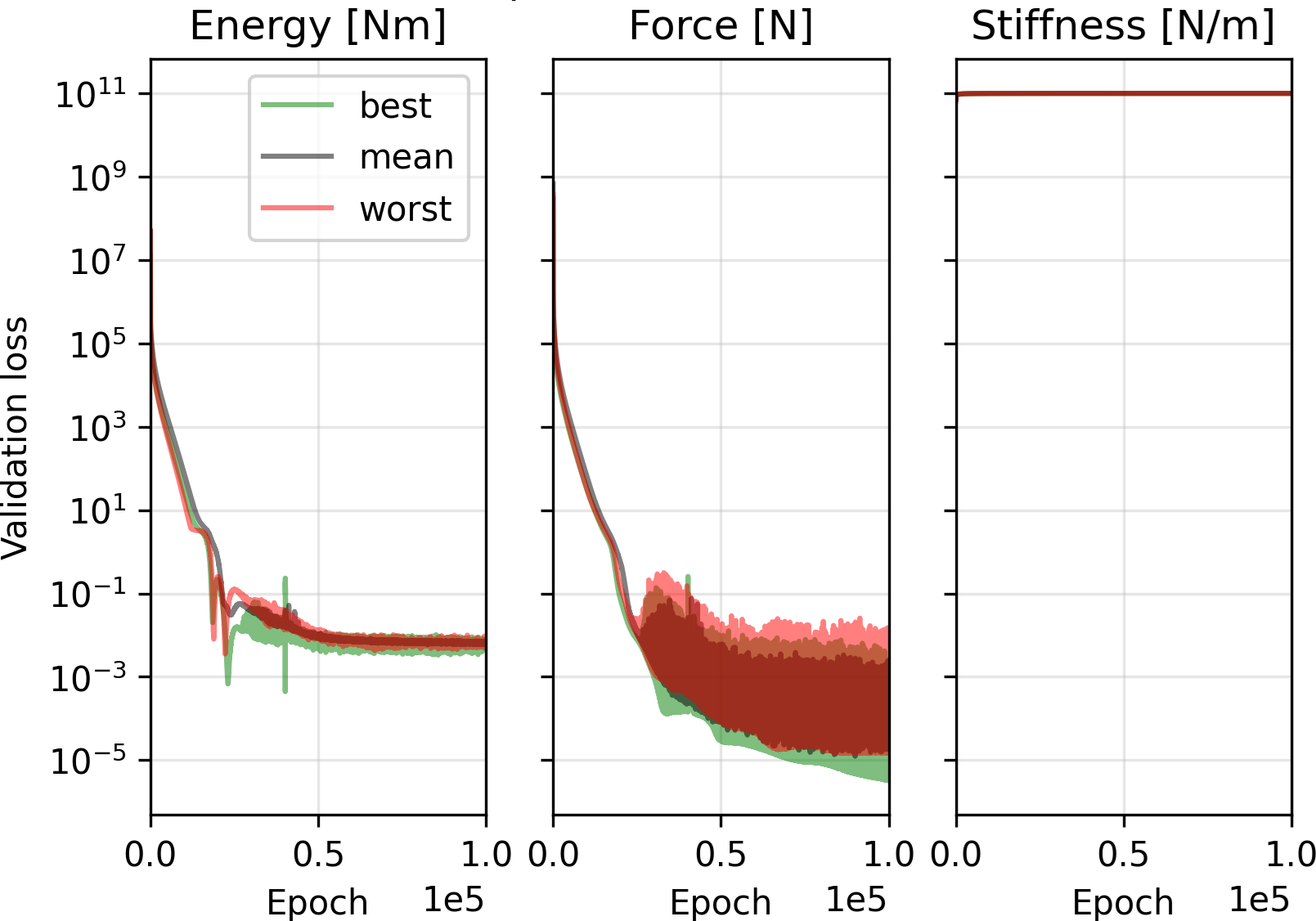}
    \caption{Validation loss curves for the chosen architecture of two hidden layers with $75$ neurons each, a learning rate of \num{1e-2}, standardized inputs, and without initializing the quadratic layer as a tangent stiffness matrix. The model is trained only using the force targets, i.e., losses for energy and tangent stiffness are just monitored.}
    \label{fig:loss_2x75_lr1e-2}
\end{figure}

In preparation for Sobolev training in the following section, \Cref{fig:loss_2x75_lr1e-2} visualizes the chosen variant's validation loss histories.
Instead of all $10$ initializations, the best, the worst, and the mean are shown.
The training exclusively uses the first derivative, but losses for energy and stiffness are additionally monitored.
It can be seen that the strain energy is learned well even though it is not part of training.
The reason is that if two functions have the same gradient, the function only differs by a constant offset.
However, this offset is set exactly by the offset for consistent energy, i.e., zero energy at zero displacement.
In contrast, the second derivative proves to be rather inert and remains constant during training.
This behavior is shared by all variants.
The loss history is rather noisy but keeps improving.
Improvements often occur in sudden steps, indicating escape from a local minimum.
The level of noise can likely be attributed to the high learning rate, which, however, leads to the best performance.

On a final note, we want to indicate two potential problems of the basic architecture used in previous studies~\cite{Asad2022,Fleres2025} and this work.
Firstly, the combination of deep architectures with unbounded, effectively polynomial activation functions leads to high initial losses.
Furthermore, a hidden layer's outputs are weighted with non-negative weights, forming the input for the next hidden layer.
This amplification grows the more hidden layers are deployed~\cite{Hoedt2023}.
In summary, these aspects limit the number of hidden layers, although tailored initializations etc.\ might counteract.
Secondly, the implementation of the non-negative constraint via Softplus can result in negligible effective weights.
Although high activations might mitigate, many connections might be inactive and can be pruned.
Alternatively, other functions for constraining parameters to be non-negative such as absolute or squaring pose aspects for future research.

\clearpage

\subsubsection{Study II: loss balancing in Sobolev training}
\label{sssec:ANN_results_Sobolev}

This second part extends the first part by Sobolev training with special emphasis on loss balancing.
Second-order Sobolev training includes targets for a function and its first two derivatives, aiming for improved performance.
Physically, these quantities correspond to strain energy, restoring forces, and tangent stiffness matrix.
All these quantities are available even from commercial \ac{FEM}.
As elaborated in \Cref{sssec:loss}, an important question is how to combine the three individual \ac{hMSE} losses.
Potential problems include differences in orders of magnitude, as demonstrated in \Cref{fig:reduced_data}, and in convergence rates.

In detail, the optimal \ac{ANN} variant determined in the previous part is augmented with the approaches for multiple losses described in \Cref{sssec:loss}.
These techniques fall in one of three categories.
The first category features simple approaches to establish a baseline.
Specific approaches include the sum of individual losses and learning the Hessian only.
The latter uses the fact that if the Hessians match, the two offsets for consistency ensure matching gradient and function.
The second category deploys a weighted sum of individual losses.
The weights aim to compensate the differences in orders of magnitude of function, gradient, and Hessian.
Three different weighting schemes are deployed, including an intuitive scheme derived from \Cref{fig:reduced_data}, a scheme based on the square ratio of maximum values~\cite{Pukrittayakamee2011}, and a scheme based on the squared standard deviation.
For all schemes, the weights are rounded to a single significant digit.
All three weighting schemes are studied in full factorial combinations with a scheduled Hessian loss and gradient aggregation via \ac{JD}.
These options and their details are provided by \Cref{tab:Sobolev_settings_weighted_sum}. 
The third category of loss balancing techniques utilizes dynamic weighting.
The goal is to additionally compensate the different losses' varying convergence rates.
Two dynamic loss balancing schemes are studied: loss balancing for \acp{PINN} and a custom scheme for same-scale losses.
Similarly to the weighted sums, both schemes are studied in all combinations with a scheduled Hessian loss and gradient aggregation via \ac{JD}.
Each of the resulting $22$ variants is initialized $10$ times.
All remaining settings remain the same as in the previous study in \Cref{sssec:ANN_results_base}.

\begin{table}[h]
\centering
\begin{threeparttable}  
\caption{Options for loss balancing based on a weighted sum of individual losses in second-order Sobolev training.}
\label{tab:Sobolev_settings_weighted_sum}
\begin{tabularx}{\textwidth}{ll}
\toprule 
Option & Values \\
\midrule 
Weighting scheme &  Intuitive ($1$, $1$, $\num{1e-9}$)\\
& Maximum ($1$, $\num{9e2}$, $\num{9e-8}$) \\
& Standard deviation ($\num{3e-1}$, $\num{2e-2}$, $\num{6e-12}$)\\
\midrule 
Schedule & None \\
& scheduled \\
& (Hessian linear ramp from \num{50000} to \num{100000} epochs)\\
\midrule 
\ac{JD} & Yes \\
& No \\
\bottomrule 
\end{tabularx}
\end{threeparttable}
\end{table}

The performance of all variants is statistically assessed analogously to the first part.
In detail, each initialization is represented by its best validation loss and its final test loss.
These data are visualized as box plots for both validation and test data.
The results and visualizations are divided into two parts.
The first part presents the first two categories of loss balancing, i.e., simple approaches and weighted sums.
The second part provides the results of dynamic loss balancing.
Note that while all three individual losses contribute to the total loss and influence the learning process, only the force or gradient loss is considered to report results. 
This decision ensures comparability with the results of the first study in \Cref{sssec:ANN_results_base}.
Moreover, it emphasizes the force as the main quantity of interest.

\begin{figure}[!htb]
    \centering
    \includegraphics[width=1\linewidth]{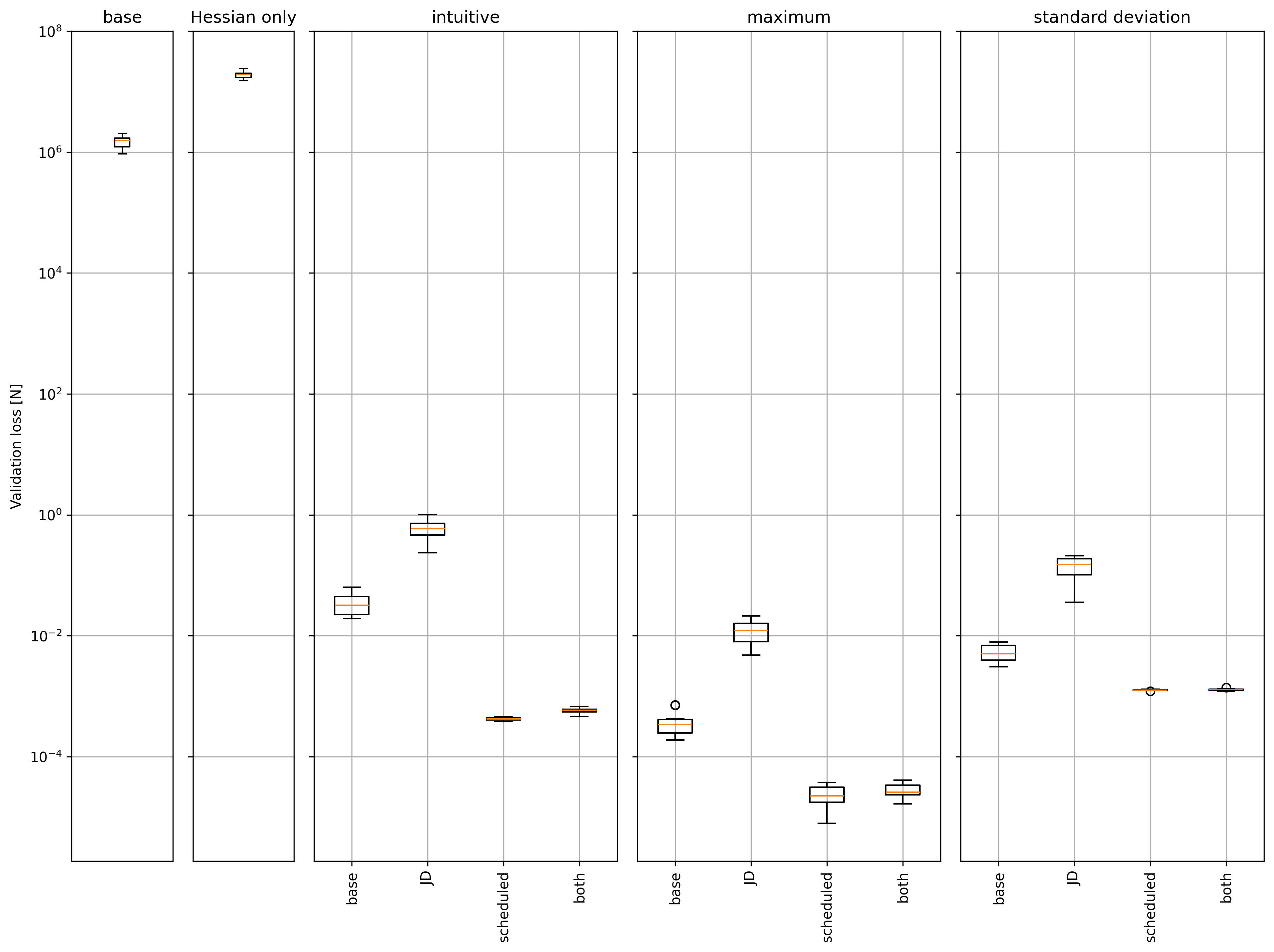}
    \caption{Statistical evaluation of model performances with respect to the validation error: each box summarizes the $10$ initializations per variant. The five subplots correspond to different static weighting schemes. While the first two are for establishing a base line only, the remaining three schemes are combined with \ac{JD}, scheduled Hessian loss, or both.}
    \label{fig:Plot_ANN_performance_Sobolev_validation}
\end{figure}

\Cref{fig:Plot_ANN_performance_Sobolev_validation} visualizes the performance of Sobolev training on the validation data using simple or weighted sum-based loss balancing.
Both simple approaches, a plain sum and Hessian only, perform about four orders of magnitude worse than weighted sums.
These approaches do not account for different orders of magnitude of the three targets.
Therefore, the stiffness loss dominates the total loss.
In addition, the hyperparameters have been optimized in \Cref{sssec:ANN_results_base} for the force loss.
Exploring the full potential of all variants would require individual studies of hyperparameters.
However, the number of variants grows exponentially.
Therefore, these investigations are beyond the scope of this work.
In contrast, the static weighting schemes perform significantly better.
Specifically, the weighting scheme based on maximum values performs best, often by more than an order of magnitude.
Its accuracy is superior in every combination of loss scheduling and \ac{JD}.
These two options show a consistent effect for all weighting schemes with loss scheduling improving performance and \ac{JD} worsening it.
The latter effect is more pronounced when \ac{JD} is deployed without loss scheduling.
We attribute the positive effect of loss scheduling to more training focused on the force loss and the fast-converging energy loss.
As the force loss constitutes the evaluation metric, the focus leads to better reported accuracy.
Similarly, \ac{JD} modifies the individual parameter updates to achieve an update that is optimal for all objectives.
Consequently, the aggregated update likely sacrifices individual improvements.
Although all losses aim to match the same reference function, their optimization directions apparently do not align.
In summary, a weighted sum based on maximum values combined with a scheduled Hessian loss achieves a superior performance.
This variant's best initialization reaches a validation loss of \qty{7.92e-06}{\newton} and a test loss of \qty{1.60e2}{\newton}.
In contrast, the best variant of part one without Sobolev training achieves a validation loss of \qty{3.27e-06}{\newton} and a test loss of \qty{2.47e2}{\newton}.
As model selection exclusively relies on the validation loss, the performance of part one is not matched by Sobolev training based on statically weighted sums.

\begin{figure}[!htb]
    \centering
    \includegraphics[width=1\linewidth]{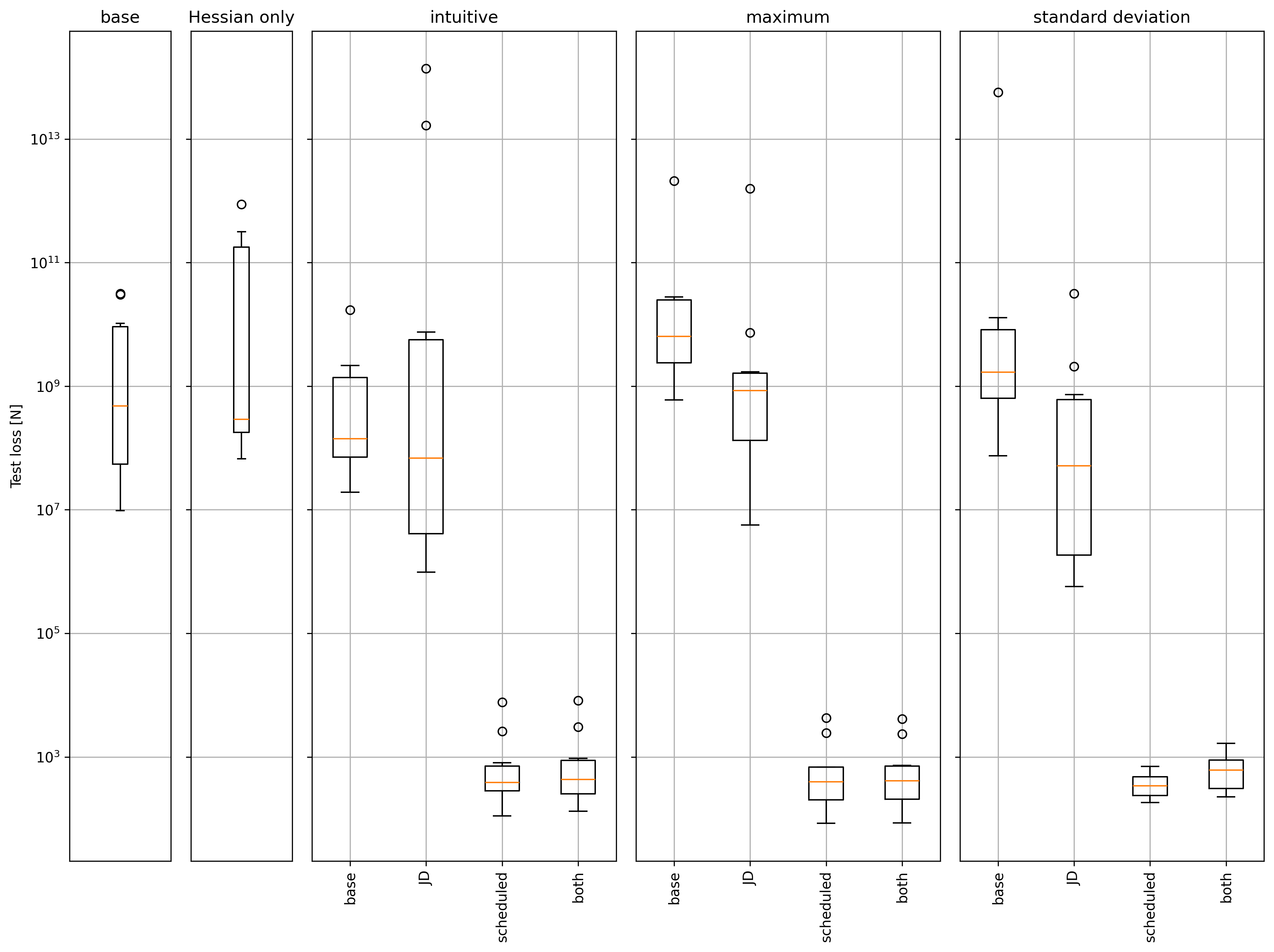}
    \caption{Analogous plot to \Cref{fig:Plot_ANN_performance_Sobolev_validation} to statistically summarize variants, but for the test error instead of validation error. Note that the test data exclusively comprise extrapolation load cases.}
    \label{fig:Plot_ANN_performance_Sobolev_test}
\end{figure}

\Cref{fig:Plot_ANN_performance_Sobolev_test} reports the variants' performances when evaluated on the challenging test data.
In general, the test losses exceed the corresponding validation losses by several orders of magnitude.
This is a result of the challenging, extrapolation-focused test set.
The performance of a simple sum and learning the Hessian become more comparable to weighted sums but still do not reach their level.
Again, loss scheduling achieves the largest improvement for each weighting scheme.
\ac{JD} improves the median and best performance if deployed on its own, but slightly worsens results if combined with loss scheduling.
The three weighting schemes vary in their performance and ranking depending on the options combined with them.
The overall best test performance is achieved by loss scheduling and either the intuitive or maximum value-based weighting scheme.

An alternative loss-balancing technique is to use dynamic weights.
In contrast to static weights that account for different scales in values, dynamic weights can potentially also account for different convergence speeds.
Two dynamic weighting schemes are studied, again in all combinations with loss scheduling and \ac{JD}.

\begin{figure}[!htb]
    \centering
    \includegraphics[width=1\linewidth]{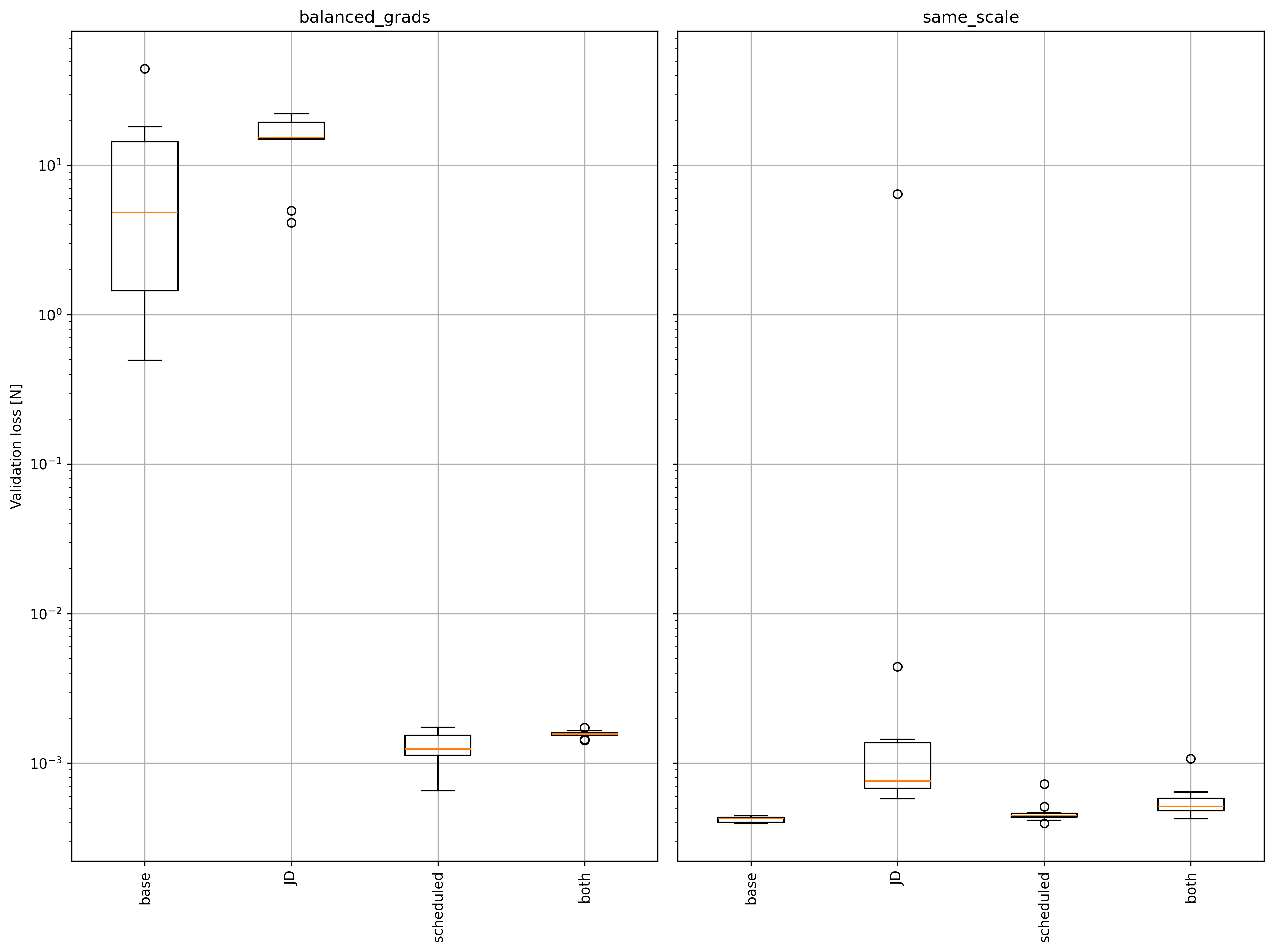}
    \caption{Statistical evaluation of model performances with respect to the validation error: each box summarizes the $10$ initializations per variant. The two subplots correspond to different dynamic weighting schemes, which are combined with \ac{JD}, scheduled Hessian loss, or both.}
    \label{fig:Plot_ANN_performance_Sobolev_validation2}
\end{figure}

\Cref{fig:Plot_ANN_performance_Sobolev_validation2} visualizes the accuracy of variants using dynamic weighting schemes on the validation data .
These variants include two weighting schemes, \ac{PINN}-inspired and custom loss balancing.
Both schemes are studied in full-factorial combinations with loss scheduling and \ac{JD}.
The custom loss scheme achieves superior accuracy for every combination.
However, a tuned learning rate might mitigate as the former amplifies low losses while the latter scales down large terms.
The effects of loss scheduling and \ac{JD} are similar as for the simple and statically weighted variants.
\ac{JD} worsens performance, both when applied on its own or combination with loss scheduling.
Loss scheduling brings major improvements for the \ac{PINN}-inspired scheme and minor degradation for the custom scheme. 
However, even the best-performing variant fails to achieve validation losses below \qty{1e-5}{\newton}, unlike the best variants trained without Sobolev training or with static weighting.

\begin{figure}[!htb]
    \centering
    \includegraphics[width=1\linewidth]{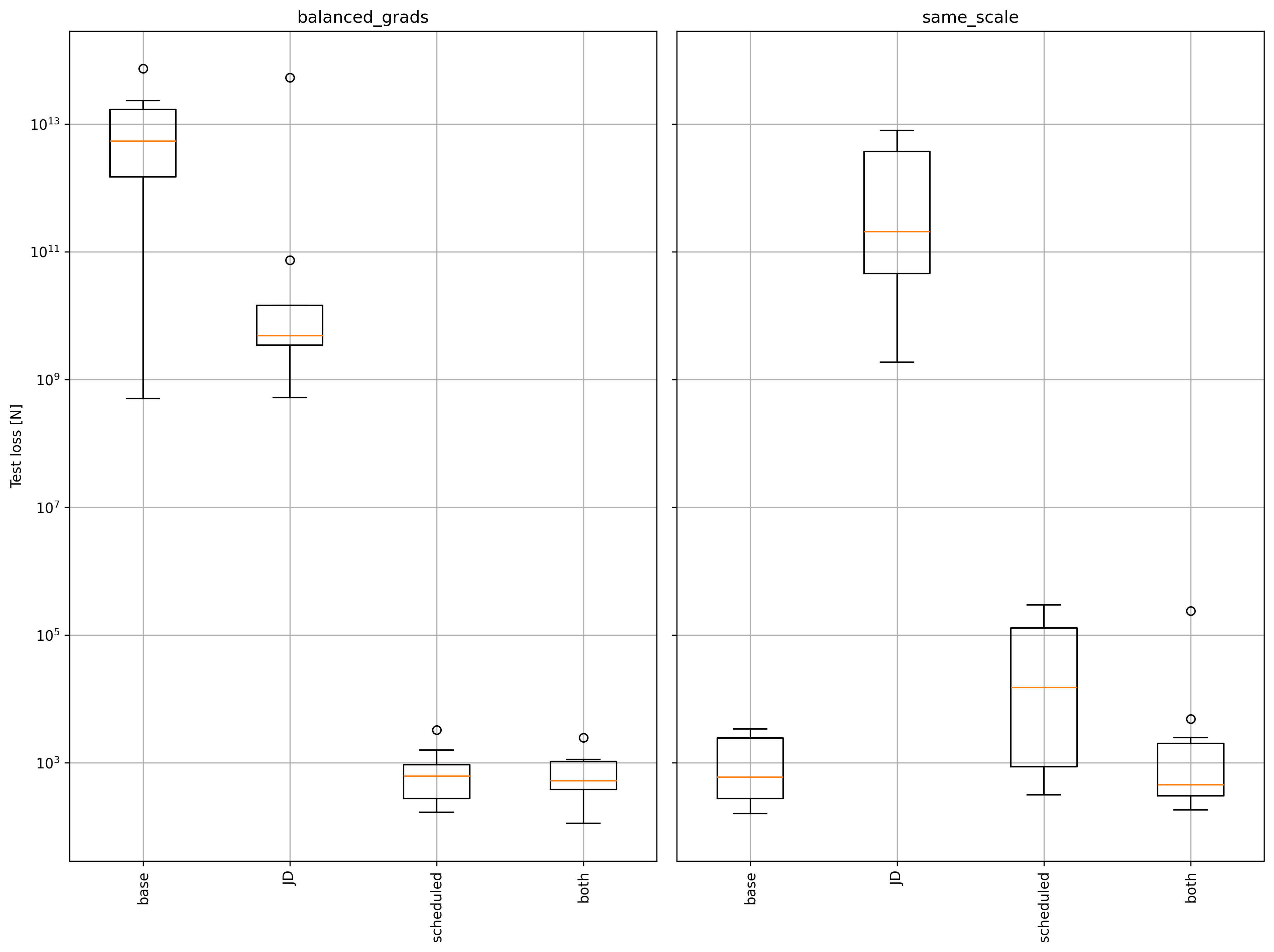}
    \caption{Analogous plot to \Cref{fig:Plot_ANN_performance_Sobolev_validation2} to statistically summarize variants, but for the test error instead of validation error. Note that the test data exclusively comprises extrapolation load cases, forming a challenging task for data-driven models.}
    \label{fig:Plot_ANN_performance_Sobolev_test2}
\end{figure}

\Cref{fig:Plot_ANN_performance_Sobolev_test2} presents the corresponding test losses.
Here, the effect of the weighting scheme, loss scheduling, and \ac{JD} is more individual.
The \ac{PINN}-inspired weighting scheme benefits from both, loss scheduling and \ac{JD}.
While the former significantly reduces spread and improves accuracy, the latter achieves superior outliers.
In contrast, the custom scheme performs well without modifications.
\ac{JD} on its own significantly worsens accuracy and also loss scheduling has a negative impact.
In combination, they outperform the unmodified scheme by a minimal margin.

In summary, the best performing \ac{PANN} is the $2 \times 75$-architecture with a learning rate of \num{1e-2}, standardized inputs, and a random initialization of the quadratic pass-through layer.
This variant's best initialization is deployed for a final evaluation in \Cref{ssec:Hyperreduction_results}.
Despite diverse approaches, second-order Sobolev-training has not improved accuracy for the selected metrics.
A possible reason is the model selection based on the loss for force prediction.
Potentially, the inert loss of the Hessian may cause tradeoffs that sacrifice accuracy on other objectives.
The chosen metrics only penalize the sacrifices, but do not reward improvements in other objectives.
However, the force is the main quantity of interest for hyperreduced models.
These application-induced priorities must reflect in model selection.
Nevertheless, the notable success of gradient-enhanced methods in other research fields suggests untapped potential beyond what has been achieved in this work.

\clearpage

\subsection{Results hyperreduction}
\label{ssec:Hyperreduction_results}

This subsection assesses the prediction accuracy of the two hyperreduced models and compares them with the reference work~\cite{Fleres2025}.
The two models deploy either \ac{TPWL} or the best performing \ac{PANN} of \Cref{ssec:PANN_results} for hyperreduction.
A \ac{PANN} of the same type was also studied in the reference work for the same numerical case study.
Therefore, the corresponding results are included to evaluate the enhancements proposed in this work.
In alignment with the reference work, the comparison considers three quantities via the  predictions' relative error magnitude.
More details, i.e., component-wise comparisons, are presented in \ref{app:Hyperreduction_results}.
Special emphasis is placed on the models’ extrapolation capabilities, as this property is of great importance for applications such as system-level simulation or feedback control.
Consequently, all results differentiate between interpolation and extrapolation, i.e., predictions for the train and validation sets and for the test set, respectively.

The comparison relies on three quantities, comprising the reduced internal force, the reduced state, and the physical displacement of the beam's tip in load direction.
These differ in their computations and partially rely on solvers that may not be consistent between studies.
The first quantity is the reduced force, which is predicted for a given reduced state.
Its computation constitutes a single function call.
Therefore, implementations are expected to hardly differ.
The second quantity is the reduced state, which is predicted for a given external load.
In contrast to the previous case, its computation poses a nonlinear system of equations that requires an iterative solution.
The corresponding solvers and their settings might deviate, likely affecting errors.
In order to mitigate this effect, this work deploys the same Newton-Raphson scheme as a solver as used in its reference, the \ac{FEM} software.
Furthermore, the analysis features the same underlying load cases as the \ac{FEM} model, using the last substep's solution as the initial guess for the next one.
The third quantity is the displacement of the beam tip, which is derived as a linear combination of the predicted reduced states.
Note that the reduced basis in the reference work also covers parametric effects compared to this work covering extended load cases.
The different reductions lead to deviations, especially for the beam tip's displacement.
However, it poses a physical quantity that is easy to interpret.

The accuracy of all predicted quantities is assessed using the relative error magnitude per sample.
For the prediction of the reduced state vector  $\vi{\hat{x}}_{r, i}$, the relative error $\varepsilon_{\vi{x}_r,i}$ is computed as
%
\begin{equation}
    \varepsilon_{\vi{x}_r, i} = \dfrac{\norm{\vi{\hat{x}}_{r, i} - \vi{x}_{r, i}}}{\norm{\vi{x}_{r, i}}} \, ,
    \label{eq:norm_error}
\end{equation}
where $i$ denotes the load or time step.
Similarly, a relative error $\varepsilon_{\vi{f}_r, i}$ is calculated for the reduced force vector, using $\vi{f}_{r, i}$ instead of $\vi{x}_{r, i}$.
In addition, the relative error for the beam's tip displacement $\varepsilon_{U_{tip}, i}$ is computed.
This unreduced quantity offers physical intuition and easier assessment of approximation quality, but also features an additional error source due to different reduced bases.
On a final note, the relative error naturally punishes deviations for small reference values.
Therefore, larger errors are expected in these regions.

Another critical aspect for comparability between models is the consistent usage of data.
The \ac{PANN} has been trained using \qty{50}{\percent} of the $101$ samples for a tip load of \qtyrange{0}{300}{\newton}.
The \ac{TPWL} model uses the same percentage of the same data set, but with equidistant sampling.
Note that for samples included in the \ac{TPWL} model's construction, interpolation reproduces the reference solution, and relative errors reach the level of numerical noise.
These samples have been excluded from the error plots, providing a conservative estimate of the \ac{TPWL} model's accuracy and clearer visualizations.
Note that the \ac{PANN} in the reference work~\cite{Fleres2025} has been trained on $80$ full trajectories of the same load case but with variations of two material parameters.
As discussed earlier, we expect that reducing the number of training samples from \num{8000} to \num{51} increases the difficulty more than the slightly simpler \ac{PANN} architecture mitigates it.

\begin{figure}[!htb]
    \centering
    \includegraphics[width=1.\linewidth]{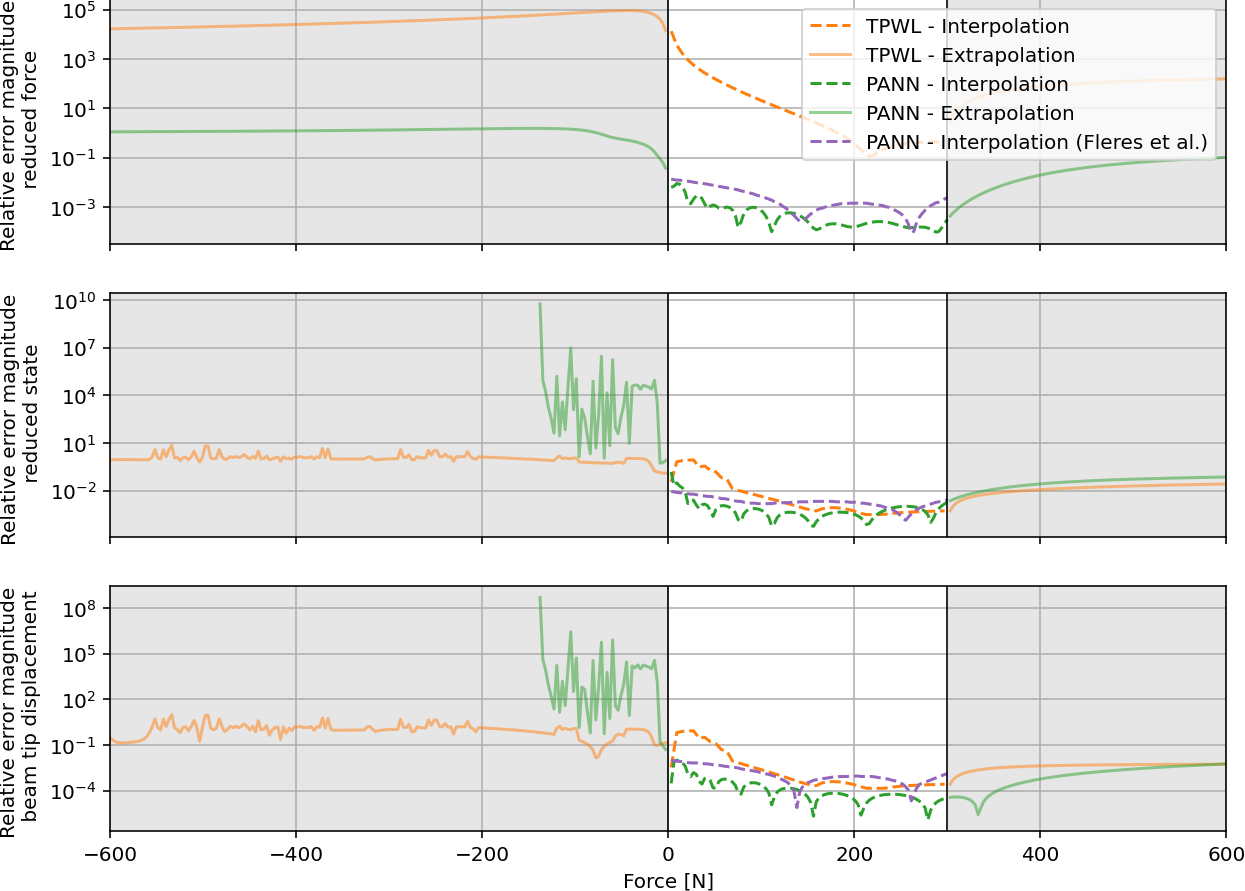}
    \caption{Relative error magnitudes of three quantities predicted by the \ac{PANN}, by \ac{TPWL}, and by the \ac{PANN} of the reference work~\cite{Fleres2025}. The quantities include the reduced force obtained by a function evaluation, the reduced state computed as an iterative nonlinear solution, and the displacement of the beam's tip, which is postprocessed from the reduced state but offers more physical intuition. Shaded regions and solid lines mark the extrapolation domains. Note that the \ac{TPWL} model includes every second sample of the interpolation region. For these samples, the error is not included as it reaches levels of numerical noise and would disturb the plot with a strong oscillating behavior. The data of Fleres et al.\ has been extracted from their Fig.~A.9~\cite{Fleres2025}. }
    \label{fig:Plot_Resultcomparison_relative_error_all}
\end{figure}

\Cref{fig:Plot_Resultcomparison_relative_error_all} presents the relative error magnitudes for the two hyperreduced models of this work and for the model from literature.
The three models are compared in predicting the reduced force, the reduced state, and the displacement of the beam's tip.
Furthermore, interpolation and extrapolation are clearly visually distinguished for precise evaluation.
The next three paragraphs focus on the three quantities.

For the reduced force, both \acp{PANN} perform around four orders of magnitude better than \ac{TPWL}.
Furthermore, the \ac{PANN} of our work generally outperforms the literature model by up to an order of magnitude.
With respect to interpolation and extrapolation, all models exhibit a similar behavior.
As expected for data-driven methods, performance for interpolation is drastically better than for extrapolation.
Extrapolation can be subdivided into two scenarios, increasing the load and loading in the opposite direction.
From these two, increasing the load of the trained load case yields smaller errors.
These deteriorate with distance from the training domain, reaching \qty{10}{\percent} for our \ac{PANN}.
The presumably related scenario of reversing the load direction causes the highest errors, quickly reaching \qty{100}{\percent} for this work's \ac{PANN}.
Note that the \ac{PANN} incorporates unbounded activation functions and non-negative weights.
Therefore, higher deviations with increasing distance to the training region would be a reasonable expectation, which, is not confirmed by experiments.
In summary, only the two \acp{PANN} achieve acceptable error levels, and even then, only for interpolation.
For extrapolation, their accuracy deteriorates rapidly.

The reduced state is the more relevant quantity, as it represents the computations of a deployed model.
The computation involves an iterative nonlinear solution where convergence is far from guaranteed for surrogate models such as plain \acp{ANN}~\cite{Fleres2025}.
Therefore, not only accuracy, but also robustness are to be compared.
In this case study, all models converge within the interpolation region.
This fact can be attributed to the mathematical properties of their tangent stiffness matrices.
The convergence continues for the extrapolation scenario of increasing the load.
However, for loading in the opposite direction, the \ac{PANN} diverges after a few substeps despite more accurate force predictions.
In this evaluation, the last substep's solution constitutes the initial guess for the next one.
Hence, deviations accumulate, potentially leading to a complete failure.
Consequently, we attribute the \ac{PANN}'s divergence to initial deviations quickly accumulating, leading the state trajectory in a wrong direction with no chance of recovery. 
In contrast, the \ac{TPWL} model exhibits extraordinary robustness.
It provides a solution for every step, even though relative errors reach \qtyrange{100}{1000}{\percent}.
The reason is likely the nearest-neighbor-like interpolation.
Even if all linearization states are far away, the \ac{TPWL} model still reverts to the closest one.
Moreover, the linearized matrix is almost always a single original sample due to sharp weighting functions, which strongly promotes the desired mathematical properties.
In terms of accuracy, the most notable change compared to the force is the drastically improved performance of \ac{TPWL}.
It not only approaches the \acp{PANN}' performances for interpolation, but even exceeds the \ac{PANN} for all extrapolations.
When comparing the two \acp{PANN}, our model almost exclusively outperforms its reference by up to an order of magnitude.
A major reason for the remarkable performance of \ac{TPWL} despite high errors in the force prediction is the original system's high stiffness.
In other words, the force is significantly more sensitive to changes in the state than vice versa.
This effect is further amplified by the stress-stiffening behavior shown in \Cref{fig:reduced_data}, i.e., more pronounced for high loads.

The last quantity is the displacement of the beam's tip.
Because it is computed as a linear combination of the reduced state, it inherits its qualitative error behavior.
Note that the two \acp{PANN} use different reduced bases, introducing an additional source of deviation and decreasing comparability.
Therefore, the comparison focuses on \ac{PANN} and \ac{TPWL}.
While the \ac{PANN}'s divergence for reverse loading remains, its extrapolation accuracy for an increasing load surpasses the \ac{TPWL} model.
This difference between errors for the reduced state and the displacement can be attributed to the varying importance of the state vector's components for computing the displacement.
Hence, while the state vector error considers all components, the displacement weights them in its computation.
Therefore, the \ac{PANN} seems to match the relevant state vector component better than the \ac{TPWL} model.

In summary, the comparison in \Cref{fig:Plot_Resultcomparison_relative_error_all} provides some interesting insights into the different models.
Both \acp{PANN} guarantee the desired mathematical properties and elegantly embed physical knowledge.
As a result, both achieve convergence for interpolation scenarios, which is far from guaranteed for \acp{ANN}~\cite{Fleres2025}.
They achieve good accuracy for interpolation with our \ac{PANN} outperforming its counterpart in most situations.
Note that our \ac{PANN} has been selected based on the experiments in \Cref{ssec:PANN_results}, which details the effect of proposed enhancements.
Another major contribution to its superior performance is the hyperparameter optimization, also emphasizing the importance of multiple initializations per variant.
In comparison, the much simpler \ac{TPWL} model achieves poor accuracy for the force, but competitive levels for the significantly more relevant state vector.
A highly important aspect where models differ significantly is convergence of state predictions in extrapolation scenarios.
This aspect is critical for typical applications of \acp{ROM}, e.g., at system-level or in control.
In both cases, the models must be able to handle load cases that deviate from those they were trained on.
A single outlying state prediction can lead the model into divergence.
Therefore, models should be robust and ideally be able to recover from disturbances.
For this case study, both \ac{TPWL} and our \ac{PANN} converge for extrapolation with a higher load.
However, the \ac{PANN} quickly diverges for loading in the opposite direction.
In contrast, \ac{TPWL} consistently converges as it always reverts to a linearization of the original model.
In summary, \ac{TPWL} achieves comparable accuracy for state vector predictions and provides unmatched robustness for a data-driven hyperreduction method.

\clearpage

\section{Conclusion}
\label{sec:Conclusion}

This work has successfully enhanced a novel non-intrusive hyperreduction method~\cite{Fleres2025} which originates from constitutive modeling~\cite{Asad2022}.
It relies on a \ac{PANN} to not promote but guarantee several desirable properties.
These properties include a positive semi-definite tangent stiffness matrix, which strongly aids convergence and stability, and consistent energy and forces.
Furthermore, it offers some physical interpretation as it learns a potential energy and derives forces and tangent stiffness matrices as derivatives.
In contrast to established hyperreduction such as \ac{DEIM} or \ac{ECSW}, the method is conveniently compatible to commercial \ac{FEM} software.
Note that the method is compatible to other physical domains, e.g., temperature-dependent thermal properties or magnetic saturation.
However, the guaranteed properties also prevent certain applications where they are not featured by target quantities.

This work discards the method's parametric capability and studies several enhancements to incorporate more physics or best practices from machine learning.
These enhancements include standardized inputs, a physical initialization of the quadratic pass-through layer, a more extensive hyperparameter study that also includes the learning rate, initializing each \ac{ANN} $10$ times to account for its stochastic nature, second-order Sobolev training, and evaluating extrapolation capabilities due to their importance for typical applications.
The idea of second-order Sobolev is to use the function, the gradient, and its Hessian in training to improve the \ac{ANN}'s performance.
In a physical sense, these quantities correspond to energy, forces, and tangent stiffness matrices.
Moreover, all are conveniently available even for commercial \ac{FEM} software.
Three target quantities lead to three losses that must be carefully balanced during training.
Also note that writing tangent stiffness matrices takes significantly longer and uses more space than energies, states, and forces.
Another extension is a comparison of the \acp{ANN} with an alternative, i.e., a model hyperreduced via \ac{TPWL}.
\ac{TPWL} is robust, requires no training, can be extended on the fly to some degree, and synergizes well with control via gain scheduling~\cite{Schuetz2025}.

All numerical experiments relied on the numerical case study of a static nonlinear cantilever beam, which was reproduced from the reference work~\cite{Fleres2025} for comparability.
The hyperparameter optimization achieved a \ac{PANN} that outperforms reference work by a factor of $25$ with respect to training loss.
Note that our \ac{PANN} discards the parametric capabilities and comprises \num{6.47e3} parameters compared to  \num{7.18e3} of the reference model.
However, its training relies on only \num{51} instead of \num{8000} samples.
Therefore, the increase in challenge dominates.
As a result, the improvement in accuracy especially highlights the necessity of multiple initializations to find different local minima.
Studying our proposed enhancements indicates beneficial effects of standardized inputs, which does not hold true for the physical initialization and second-order Sobolev-training.
Nevertheless, the success of Sobolev training for \acp{ANN} and gradient-enhanced methods in other fields indicates potential not yet reached within this work.

After studying hyperparameters and enhancements, the accuracy of two models was compared with the reference work in a final evaluation.
These models comprise the best-performing \ac{PANN} and a plain \ac{TPWL} model.
The \ac{PANN} achieves great accuracy and mostly surpasses its reference and \ac{TPWL} for interpolation despite training on only \num{51} samples.
Therefore, it qualifies for regression where samples are scarce.
Furthermore, it converged when deployed in a nonlinear solver for all interpolation tasks.
This achievement is not guaranteed for \acp{ANN} and is associated with its numerical properties.
\ac{TPWL} deviates significantly for force predictions, but reaches competitive accuracy when deployed in a hyperreduced model to compute the system's state.
The different levels of accuracy for the two quantities arise due to the mechanical stiffness, which causes deviating sensitivities.

In addition, this work studied extrapolation as it is critical for desired applications such as system-level simulation or control.
In both, models likely encounter states that differ from training.
However, the \ac{PANN} exhibits poor extrapolation capabilities and quickly diverges for presumably related load cases, preventing its application.
Therefore, we recommend \ac{TPWL} over \acp{PANN} for hyperreduction in scenarios where convergence for extrapolation is critical.
Note that extrapolation leads to rapidly deteriorating accuracy due to the data-driven nature of these methods.

In our opinion, these statements reflect the current state of the art, but further enhancements may resolve existing drawbacks.
\Cref{sec:Outlook} presents the corresponding suggestions for future work.

\clearpage

\section{Outlook}
\label{sec:Outlook}

While the \ac{PANN}-based hyperreduction elegantly guarantees desired properties and incorporates physics by design, it fails for extrapolation.
To mitigate this drawback, the remainder of this section proposes aspects to investigate.
These aspects are categorized into sampling strategies and \ac{PANN} enhancements.

\subsection{Sampling finite element models}
Sampling provides an easy solution to insufficient performance for extrapolation.
Note that a reliable operation requires anticipating all potential trajectories.
Moreover, sampling is costly as the high-dimensional \ac{FEM} must be solved.
However, these costs are decoupled from the final model's operation and can be parallelized.

A convenient approach for efficiency is to sample by prescribing states~\cite{Muravyov2003, Kim2013} instead of solving the original model.
Instead of solving a nonlinear system of equations in numerous iterations, only a function must be evaluated as the state solution is completely determined.
As a result, this approach circumvents expensive nonlinear iterations, tolerances, and convergence issues.
The states to be prescribed can be expanded random reduced state vectors~\cite{Gabbay1998}.
An initial reduced basis can be obtained via an initial \ac{POD} or by a modal reduced basis~\cite{Gabbay1998} augmented by modal derivatives~\cite{Rutzmoser2017}.
A related technique for efficient sampling are nonlinear stochastic Krylov training sets~\cite{Rutzmoser2017}, which however, construct load cases and not prescribed states.
A more general and compatible sampling procedure is to deploy a greedy scheme in combination with a gaussian process.
This approach sample at maximum uncertainty and uses available resources more efficient.

Another aspect are tangent stiffness matrices which are required for second-order Sobolev training.
The corresponding files are much larger with long export and import times.
Although these steps are associated with the offline phase, efficiency can be improved by using only partial sampling~\cite{Feng2022c}.
In other words, the expensive tangent stiffness matrix is not included for every sample but only a certain percentage.

\subsection{Physics-augmented neural network}

More aspects for future work aim to improve the underlying \ac{ANN} architecture and include parameter initialization, activation functions, and parameter constraints.

The first aspect is how to initialize the \ac{ANN}'s parameters.
Established initializations rely on certain assumptions and are often derived for certain activation functions~\cite{Glorot2010,He2015} or architecture~\cite{Cranmer2020,Hoedt2023}.
Due to the SoftplusSquared activation function and  non-negative weights, these assumptions do not hold true.
Therefore, deriving an initialization scheme by stochastically analyzing the \ac{ANN} or by data-driven discovery~\cite{Cranmer2020} could lead to faster training and better performance.

The second aspect are activation functions.
These must fulfill several conditions to ensure convexity, but still achieve enough representation power and avoid exploding gradients.
Viable alternative options to SoftplusSquared include leaky Softplus or higher powers of it.
However, the latter may tend to exploding gradients.
Hence, it potentially requires careful initializations and benefits from shallow architectures.
Another approach for activation functions is to change the parameterization.
While the current implementation features a single parameter for the whole \ac{ANN}, parameters can be introduced per layer or even per neuron.
In addition, the current parameterization covers sharpness, but not the slope.
Although this parameter could be represented by weights in the subsequent layer, it may allow for faster training.
As the loss might be very sensitive to the slope parameter, a separated learning rate for this type of parameter might be appropriate.

The third aspect is constraining parameters to be non-negative.
The current implementation uses Softplus, which might lead to dying weights that do not contribute and can be pruned.
Possible alternative implementations are clipping, the absolute function, or squaring the weights.
Especially the last two options seem promising as they are continuous, likely to remain active, and do not interfere with gradient descent.

\section*{Acknowledgments}

\subsection*{Funding:}
This work is part of the project ADOPT - Adaptive Optics for THz within the priority program KOMMMA - Cooperative Multistage Multistable Microactuator Systems which is supported by the Deutsche Forschungsgemeinschaft (German Research Foundation) [grant number 424616052].

\section*{Author contributions}

\textbf{Arwed Schütz:} Conceptualization, Methodology, Software, Investigation, Writing - Original Draft, Visualization, Funding acquisition.
\textbf{Lars Nolle:} Methodology, Writing - Review \& Editing, Supervision.
\textbf{Tamara Bechtold:} Resources, Writing - Review \& Editing, Supervision, Funding acquisition.

\section*{Declaration of competing interest}

The authors declare that they have no known competing financial interests or personal relationships that could have appeared to influence the work reported in this paper.

\section*{Acronyms}

\begin{acronym}[DIPANN]
\acro{Adam}[Adam]{adaptive moment estimation}
\acro{AI}[AI]{artificial intelligence}
\acro{ANN}[ANN]{artificial neural network}
\acro{BC}[BC]{boundary condition}
\acro{BEM}[BEM]{boundary element method}
\acro{DAE}[DAE]{differential algebraic equation}
\acro{DEIM}[DEIM]{discrete empirical interpolation method}
\acro{DIPANN}[DIPANN]{derivative-informed physics-augmented neural network}
\acro{ECSW}[ECSW]{energy conserving mesh sampling and weighting}
\acro{FE}[FE]{finite element}
\acro{FEM}[FEM]{finite element method}
\acro{GENN}[GENN]{gradient-enhanced neural network}
\acro{hMSE}[hMSE]{half mean squared error}
\acro{IC}[IC]{initial condition}
\acro{ICNN}[ICNN]{input convex neural network}
\acro{JD}[JD]{Jacobian descent}
\acro{LSTM}[LSTM]{long short-term memory}
\acro{MEMS}[MEMS]{microelectromechanical system}
\acro{MIMO}[MIMO]{multiple-input multiple-output }
\acro{ML}[ML]{machine learning}
\acro{MLP}[MLP]{multilayer perceptron}
\acro{MOR}[MOR]{model order reduction}
\acro{MSE}[MSE]{mean squared error}
\acro{ODE}[ODE]{ordinary differential equation}
\acro{PANN}[PANN]{physics-augmented neural network}
\acro{PCGrad}[PCGrad]{projecting conflicting gradients}
\acro{PDE}[PDE]{partial differential equation}
\acro{PINN}[PINN]{physics-informed neural network}
\acro{POD}[POD]{proper orthogonal decomposition}
\acro{ReLU}[ReLU]{rectified linear unit}
\acro{RBF}[RBF]{radial basis function}
\acro{RMSE}[RMSE]{root mean square error}
\acro{ROM}[ROM]{reduced order model}
\acro{RWF}[RWF]{random weight factorization}
\acro{SVD}[SVD]{singular value decomposition}
\acro{TPWL}[TPWL]{trajectory piecewise linear}
\end{acronym}

\appendix

\section{Number of Parameters}
\label{app:N_par}

The total number of parameters $N_{par}$ of a \ac{PANN} used within this work is given by
\begin{equation}
    N_{par} = \underbrace{(N_I +1) \, N_N}_{\text{Input layer}} + \underbrace{N_N+N_I + N_I^2}_{\text{Output layer}} + \underbrace{2}_{\alpha, \beta} +\underbrace{(N_N+N_I+1) \, N_N \, (N_L-1)}_{\text{Hidden layers}} \, ,
\end{equation}
where $N_I$ is the number of inputs, $N_N$ the number of neurons per hidden layer, and $N_L$ the number of hidden layers.
The assumptions include a scalar output, a single quadratic pass-through layer ($N_I^2$), two global parameters for non-negative constraints and activation function sharpness, as well as removed bias for the output due to the offset to be subtracted.

\clearpage

\section{Extended hyperreduction results}
\label{app:Hyperreduction_results}

\begin{figure}[!htb]
    \centering
    \includegraphics[width=1.\linewidth]{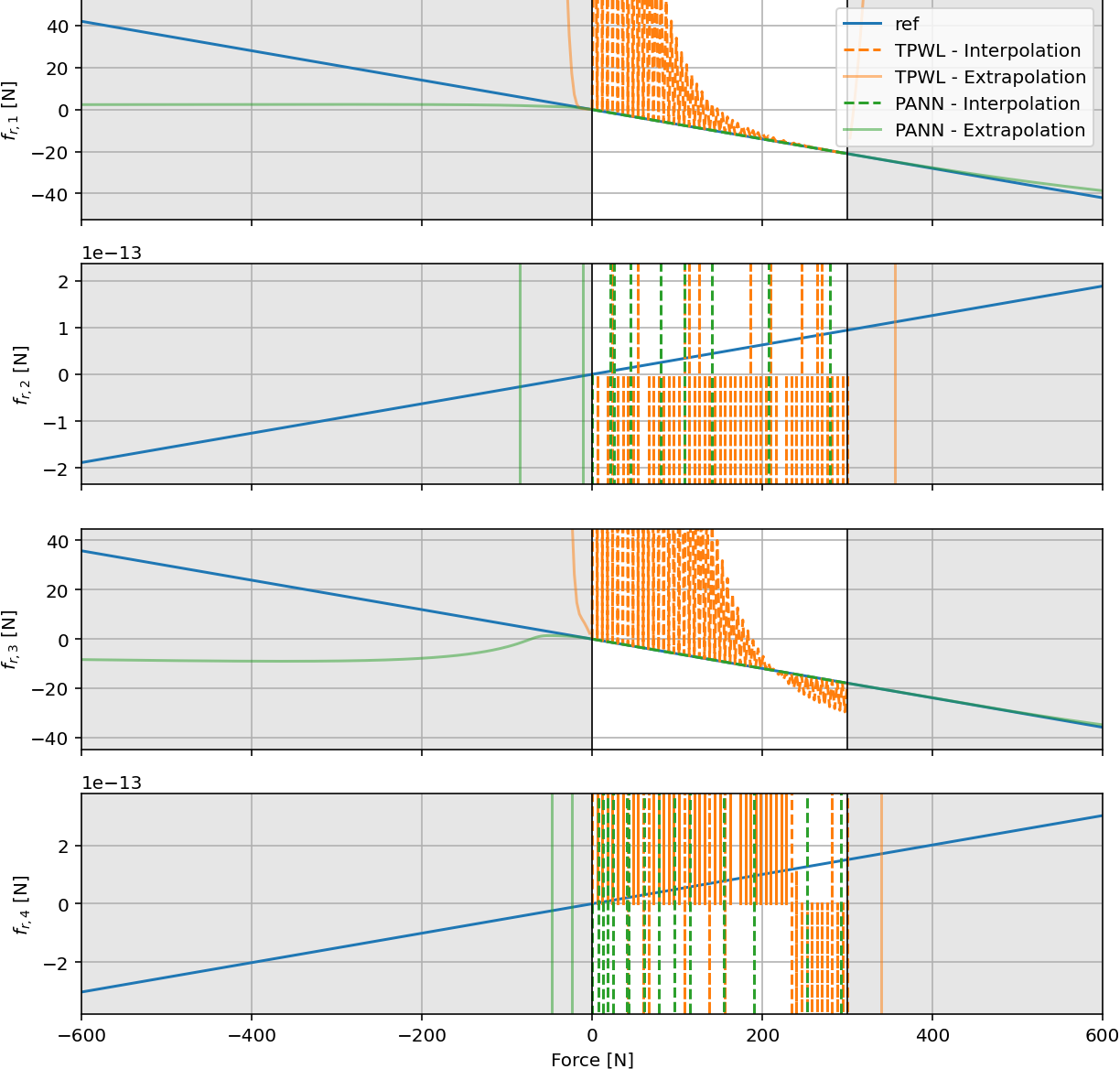}
    \caption{Reduced force components computed by the reference model, the \ac{PANN}, and by \ac{TPWL}. Shaded regions and solid lines mark the extrapolation domains. \ac{TPWL}'s oscillating behavior results from using every second sample within the training interval for its construction. Note that the second and fourth component have a magnitude of \num{1e-13} and thus indistinguishable from numerical noise.}
    \label{fig:Plot_Resultcomparison_internal_force_history}
\end{figure}


\begin{figure}[!htb]
    \centering
    \includegraphics[width=1.\linewidth]{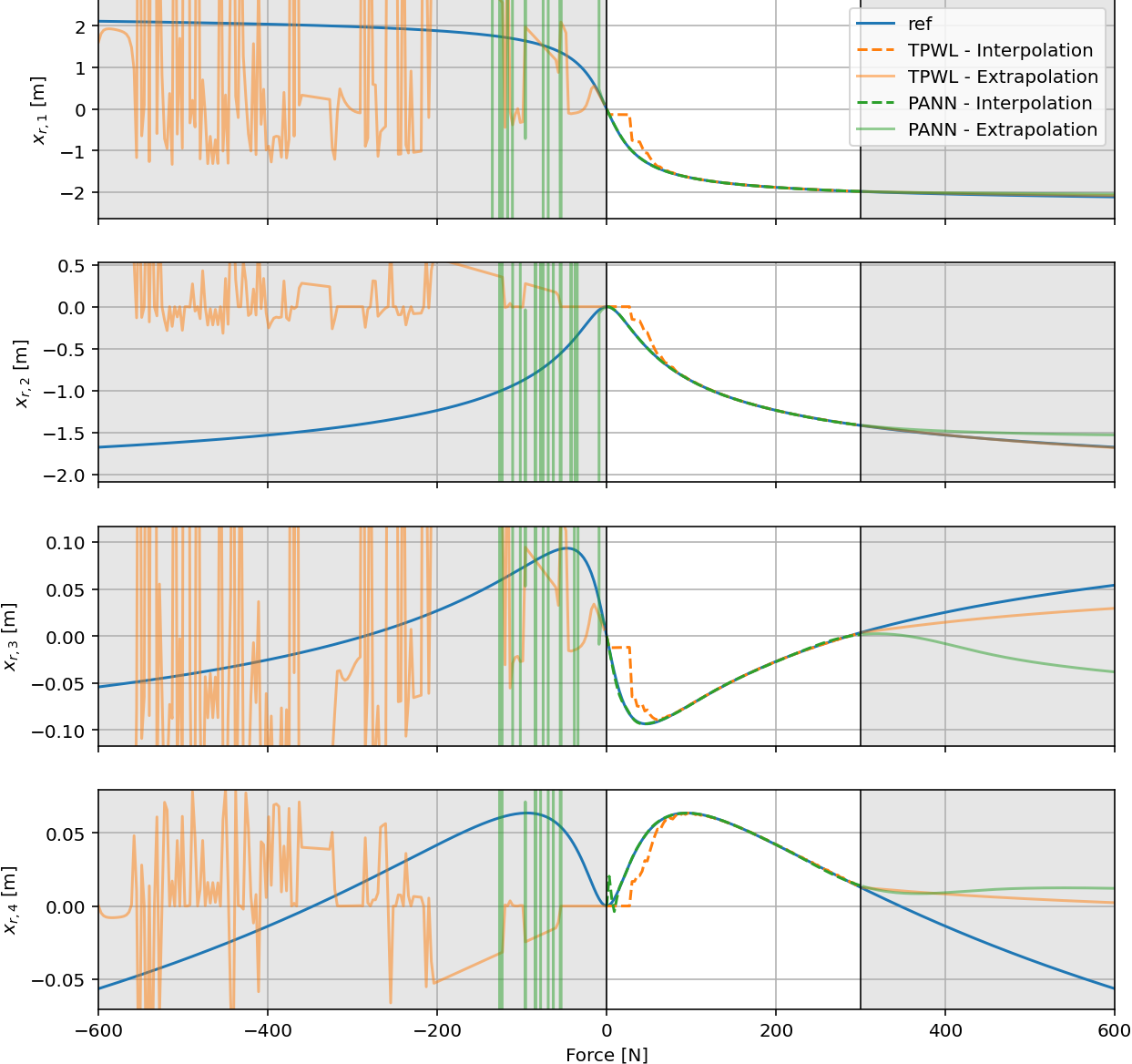}
    \caption{Reduced state components computed by the reference model, the \ac{PANN}, and by \ac{TPWL}. Shaded regions and solid lines mark the extrapolation domains. \ac{TPWL}'s oscillating behavior results from using every second sample within the training interval for its construction.}
    \label{fig:Plot_Resultcomparison_state_history}
\end{figure}


\begin{figure}[!htb]
    \centering
    \includegraphics[width=1.\linewidth]{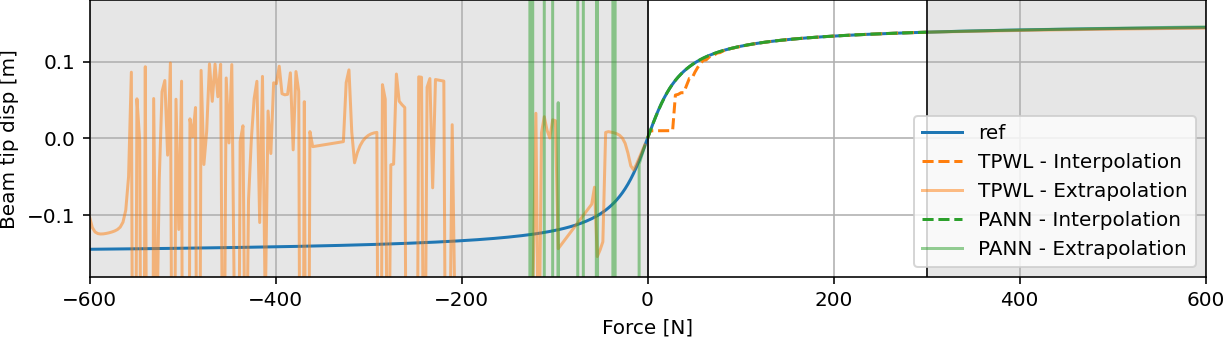}
    \caption{Beam tip displacement in load direction computed by the reference model, the \ac{PANN}, and by \ac{TPWL}. Shaded regions and solid lines mark the extrapolation domains.}
    \label{fig:Plot_Resultcomparison_output_history}
\end{figure}

\clearpage

\bibliographystyle{abbrvnat} 
\bibliography{references}

\end{document}